\documentclass[bj,numbers]{imsart}
\RequirePackage{amsthm,amsmath,amsfonts,amssymb}
\usepackage{graphicx}
\RequirePackage[colorlinks,citecolor=blue,urlcolor=blue]{hyperref}
\RequirePackage{graphicx}
\usepackage[title]{appendix}
\usepackage{mathrsfs}
\usepackage{dutchcal}
\usepackage{booktabs} 
\usepackage{caption}  
\usepackage{threeparttable} 
\usepackage{diagbox}

\usepackage{float}   
\usepackage{caption} 
\startlocaldefs
\numberwithin{equation}{section}
\newtheorem{thm}{Theorem}[section]

\numberwithin{equation}{section}

\newtheorem{lemma}[thm]{Lemma}
\newtheorem{remark}[thm]{Remark}

\newtheorem{proposition}[thm]{Proposition}
\newenvironment{Proof}{\begin{proof}}{\end{proof}}

\newcommand{\bv}{\mathbf {\widetilde v}}
\newcommand{\bV}{\mathbf {\widetilde V}}
\newcommand{\bS}{\mathbf {S}}
\newcommand{\by}{\mathbf {y}}
\newcommand{\bY}{\mathbf {Y}}
\newcommand{\be}{\mathbf {e}}
\newcommand{\bmC}{\boldsymbol{\mathbcal{C}}}

\newcommand{\diag}{\operatorname{diag}}

\newcommand{\bI}{\mathbf{I}}

\newcommand{\bbu}{\mathbf{u}}

\newcommand{\bx}{\mathbf{x}}

\newcommand{\bSigma}{\boldsymbol{\Sigma}}

\allowdisplaybreaks
\endlocaldefs

\begin{document}
	
	\begin{frontmatter}
		\title{Limiting eigen-structure of spiked sample covariance matrices under missing observations}
		\runtitle{Spiked sample covariance matrices under missing observations}
		\begin{aug}
			\author[A]{\fnms{Haotian} \snm{Cheng}
				\ead[label=e1]{chenght@stu.cqu.edu.cn}
			}	
			\author[B]{\fnms{Huiqin} \snm{Li}
				\ead[label=e2]{lihq118@nenu.edu.cn}}	
			\author[B]{\fnms{Yanqing} \snm{Yin}
				\ead[label=e3]{yinyq799@nenu.edu.cn}}	
			\author[C]{\fnms{Zhixiang} \snm{Zhang}
				\ead[label=e4]{zhixzhang@um.edu.mo}}

			\runauthor{Cheng, Li, Yin and Zhang}
			
			\address[A]{School of Mathematics and Statistics,
				Chongqing University \printead[presep={,\ }]{e1}}\address[B]{School of Statistics and Data Science,
				Nanjing Audit University\printead[presep={,\ }]{e2,e3}}
			\address[C]{Department of Mathematics,
				University of Macau\printead[presep={,\ }]{e4}}

		\end{aug}
		
		\begin{abstract}
			High-dimensional Principal Component Analysis (PCA) has become an essential tool in modern data analysis, offering dimensionality reduction and feature extraction. However, the presence of missing data introduces significant challenges, distorting the performance of PCA and complicating statistical inference. In this paper, we study the asymptotic behavior of PCA under a spiked population model with missing observations, leveraging recent advances in random matrix theory.
			We demonstrate that while the spiked sample eigenvalues exhibit asymptotic normality, the limiting parameters differ substantially from those in the complete data case, reflecting the non-trivial influence of the missing data mechanism. As an application of our results, we propose a test to evaluate the independent structure of a spiked population.
		\end{abstract}
		
		
		
	\end{frontmatter}
	
	\section{Introduction}
	With the advent of the big data era, high-dimensional data has become fundamental to various fields, including machine learning, image processing, genomics, and finance. PCA is a widely used method for analyzing such data, offering powerful tools for dimensionality reduction and feature extraction. However, high-dimensional PCA encounters significant challenges in the presence of missing data, a pervasive issue in real-world applications. Missing data can arise for various reasons: image data may be lost during transmission, financial data might become unavailable under extreme market conditions, or system incompatibilities during mergers and acquisitions may result in data loss during integration.
	
	To address missing data, two common approaches are typically employed in statistics: (1) discarding samples with missing values and (2) imputing missing values using relevant variables, followed by conducting statistical analysis on the completed dataset. However, when a substantial proportion of samples are incomplete, the first approach becomes impractical, while imputation methods often fail to fully preserve the underlying structure of the data. These challenges are particularly acute in high-dimensional settings, where PCA is frequently applied. Missing data can distort the performance of PCA, necessitating a deeper understanding of its behavior when applied to incomplete datasets.
	
	PCA relies on the eigenvalues and eigenvectors of the sample covariance matrix, which are critical tools for multivariate statistical inference \cite{1958An}. Yet, in high-dimensional settings, the eigenvalues and eigenvectors of the sample covariance matrix often diverge substantially from their counterparts in the population covariance matrix, thereby undermining PCA’s reliability. This limitation underscores the importance of studying the asymptotic behavior of these quantities in high-dimensional contexts to enable robust statistical inference when conventional methods are insufficient.
	
	The spiked population model provides a valuable framework for addressing these issues. It focuses on population covariance matrices with extreme eigenvalues and captures a key phenomenon: when certain eigenvalues of the population covariance matrix exceed a critical threshold, the corresponding eigenvalues of the sample covariance matrix deviate significantly from its bulk spectrum. These outlying eigenvalues, known as spiked eigenvalues, encapsulate crucial information about the underlying data structure. The spiked model has found widespread applications across diverse fields such as ecology, epidemiology,  finance, and climate science, where it is used to study extreme environments, high-incidence disease populations, extreme financial losses, and severe weather events.
	
	The development of the spiked model in statistical research traces back to the seminal work of Johnstone  \cite{johnstone2001distribution}, who introduced a specialized spiked model characterized by a  $p \times p$  diagonal population covariance matrix:
	$$\boldsymbol{\Sigma} = \text{diag}\{ \alpha_1, \dots,  \alpha_N, 1, \dots, 1\},$$
	where  $N$  is a fixed number representing the quantity of spiked eigenvalues. Building on Johnstone’s work, Baik and Silverstein \cite{2006Eigenvalues} made a significant contribution by demonstrating that the largest eigenvalue of the sample covariance matrix converges almost surely to a limiting value determined by boundary conditions. Bai and Yao \cite{bai2008central, bai2012sample} extended the spiked model by studying a block-diagonal population covariance matrix
	\begin{equation*}
		\boldsymbol{\Sigma} =
		\begin{pmatrix}
			\boldsymbol{\Sigma}_N & \mathbf{0} \\
			\mathbf{0} & \boldsymbol{\Sigma}_S
		\end{pmatrix},
	\end{equation*}
	and presented the central limit theorem (CLT) for {distant} spiked eigenvalues.   {In the generalized spiked model they proposed, they classify the spikes into distant spikes and closed spikes based on the sign of 
		$\psi'(\alpha_{k})$, where $\psi(\alpha_k)=\alpha_k\left[1+c\int t(\alpha_k-t)^{-1}H(dt)\right]$, $H$ is the limit spectral distribution of $\boldsymbol{\Sigma}$. If $\psi'(\alpha_k)>0$ then the corresponding sample spiked eigenvalues have limits outside the support of 
		$F^{c,H}$, the limit of the empirical spectral distribution of the sample covariance matrix $\bS_n$. In this case, they refer to these spikes as distant spikes, and in this paper, we term the corresponding eigenvectors distant spiked eigenvectors.} Then, \cite{2018Generalized} contributed by removing the block-diagonal covariance matrix assumption of the spiked model and replacing the finite fourth moment condition with a tail probability condition. Zhang et al. \cite{zhang2022asymptotic} further extended these results by considering a general form of the population covariance matrix under the finite fourth moment assumption, deriving the joint asymptotic distribution of the sample spiked eigenvalues and the linear spectral statistics. 
	
	Alongside these developments, research on spiked eigenvectors also gained prominence. Paul was the first to investigate the asymptotic behavior of spiked sample covariance matrix eigenvectors under Gaussian data conditions \cite{2007Asymptotics}. Following this, many papers have explored the properties of spiked eigenvectors under specific models. For instance, by assuming a block-diagonal structure, \cite{Morales-Jimenez2021} studied the spiked correlation matrix model, while \cite{Bao2022} examined the fluctuation of eigenvectors of the spiked covariance matrix when the bulk component is the identity matrix.
	These contributions collectively advanced the understanding of spiked models, both in terms of eigenvalues and eigenvectors,  significantly impacting various fields where high-dimensional statistics are applied.
	
	While significant progress has been made in understanding high-dimensional PCA, the impact of missing data remains largely unexplored. Drawing on recent advancements by \cite{li2024spectral}, who studied the spectral properties of random matrices with missing data using random matrix theory, this paper takes a novel approach by examining the limiting eigen-structure of the sample covariance matrix under a spiked model with missing observations.
	
	Our primary goal is to investigate the limiting behavior of distant spiked eigenvalues and eigenvectors in the presence of missing data. Specifically, we analyze how the missing data mechanism influences the eigen-structure of the sample covariance matrix. Under suitable conditions, we demonstrate that the distant spiked eigenvalues converge to a normal distribution, suggesting that even with missing data, the asymptotic behavior of eigenvalues remains well-defined.
	In addition, we establish a limiting theorem for the eigenvectors, showing that their asymptotic behavior also follows a predictable pattern in the missing spiked population model. However, a key finding is that, while the eigenvalues and the inner product of eigenvectors exhibit asymptotic normality, the limiting parameters differ significantly from those in the complete data case, revealing a non-trivial influence of the missing data mechanism.
	
	These results not only deepen the theoretical understanding of spiked models in the context of missing data but also provide a valuable foundation for statistical inference in practical applications involving high-dimensional datasets. By addressing this gap, we aim to inspire further research and new methodologies in high-dimensional statistical analysis. {As an application of our results, we propose a test to evaluate the independent structure of a spiked population.}
	
	The organization of the remaining sections is as follows: Section \ref{sec2} introduces the model assumptions and the main theorems of this paper. Simulation results for the these theorems are reported in Section \ref{se3}. Section \ref{sec3} presents an application. The proof of the main theorems are detailed in Section \ref{sec4}. Finally, the lemmas used in the analysis are provided in Section \ref{sec5}.
	
	Throughout this paper, we use the following notations: $\mathbf{A}^T$ represents the transpose of the matrix $\mathbf{A}$, $\|\mathbf{A}\|$ denotes the spectral norm of $\mathbf{A}$, $\mathbf{A} \otimes \mathbf{B}$ represents the Kronecker product of matrices $\mathbf{A}$ and $\mathbf{B}$, and $\mathbf{A} \circ \mathbf{B}$ represents the Hadamard product of matrices $\mathbf{A}$ and $\mathbf{B}$. The notation $\rm{diag}\left\{\alpha_1,\dots, \alpha_p\right\}$ refers to a diagonal matrix with $\alpha_1, \cdots,  \alpha_p$ as its diagonal elements, and $\boldsymbol{1}_n$ denotes an $n$-dimensional column vector with all elements equal to 1.   {For a matrix $\mathbf{A}=(a_{ij})_{N\times N}$, define {\rm vech}($\mathbf{A}$)=$(a_{11},\cdots,a_{N1},a_{22},\cdots,a_{N2},\cdots,a_{NN})^T$. For a sequence
		$(a_n)_{n\geq1}$ of scalars, we write $a_n=O_{p}(1)$ if $(a_n)_{n\geq1}$ is bounded in probability and  $a_n=o_{p}(1)$
		if $(a_n)_{n\geq1}$ converges to zero in probability. Similarly, we denote \(a_n = O_{\text{a.s.}}(1)\) when \((a_n)_{n\geqslant1}\) is almost surely bounded and \(a_n = o_{\text{a.s.}}(1)\) when it converges to zero almost surely.}
	\section{Model assumptions and main results}\label{sec2}
	In this section, we introduce the spiked population model under missing observations and present the main theorems of this paper. Building on the framework established in \cite{bai2012sample, Morales-Jimenez2021}, we adopt a block-diagonal spiked structure to model the population covariance matrix. To account for the random pattern of missing data, we modify the original population data by introducing a dimension-wise random factor. Specifically, each dimension is scaled by an independent Bernoulli random variable, effectively modeling the presence or absence of data in that dimension. This adjustment allows us to rigorously analyze the spectral properties of the sample covariance matrix under a realistic missing data mechanism. Specifically, the population under random missing data is formed as follows: $$\mathbf{y}=\mathbcal{b}\circ\mathbf{x}=\left(\mathbf{y}_N^T ,\mathbf{y}_S^T\right)^T,$$
	where $\mathbcal{b} =\left({b}_1,\dots,{b}_p \right)^T
	=
	\left(\mathbcal{b}_N^T,\mathbcal{b}_S^T\right)^T$ with  {  ${b}_i\sim {\rm Ber} (\theta_i)$, } $i=1, \ldots, p$, and $N$ is a fixed number. 
	Here  $\mathbf{x}=\left(\mathbf{x}_N^T ,\mathbf{x}_S^T\right)^T$ is a $p$-dimensional random vector with  mean zero and covariance matrix
	\begin{equation*}
	\boldsymbol{\Sigma} =
	\begin{pmatrix}
		\boldsymbol{\Sigma}_N & \mathbf{0} \\
		\mathbf{0} & \boldsymbol{\Sigma}_S
	\end{pmatrix}.
	\end{equation*}
	where $\boldsymbol{\Sigma}_N$ and $\boldsymbol{\Sigma}_S$ are non-random, positive semi-definite matrices with dimensions $N\times N$ and $S\times S$, respectively. Additionally, we assume that $\mathbf{x}_N$ and $\mathbf{x}_S$ are independent. Moreover, $  {\bx}$ obeys the independent component structure, that is to say,   {$\mathbf{x}=\boldsymbol{\Sigma}^{1/2}\mathbcal{x}$} with the entries of $\mathbcal{x}$ being independent and identically distributed standard random variables.   {This is a classical assumption in the field of random matrix theory, and numerous classical results, including the Marchenko-Pastur (M-P) law, are built upon the assumption that the data have an independent component structure.}  
	It is apparent that $\mathbb{E}(\mathbf{y})=\mathbf{0}$ and the covariance matrix of $\mathbf{y}$ is 
	$$
	\widetilde{\boldsymbol{\Sigma}} =
	\begin{pmatrix}
		\widetilde{\boldsymbol{\Sigma}}_N & \mathbf{0} \\
		\mathbf{0} & \widetilde{\boldsymbol{\Sigma}}_S
	\end{pmatrix}
	=\mathbcal{P}\boldsymbol{\Sigma}\mathbcal{P}+(\mathbcal{P}-\mathbcal{P}^2)\circ\boldsymbol{\Sigma}, \ {\rm{where}}	\
	\mathbcal{P} =
	\rm{diag}\left\{\theta_1,\dots,\theta_p\right\}
	=
	\begin{pmatrix}
		\mathbcal{P}_N & \mathbf{0} \\
		\mathbf{0} & \mathbcal{P}_S
	\end{pmatrix}. 
	$$
	Assume that ${\boldsymbol{\Sigma}}_N$ has eigenvalues
	$ \alpha_1, \alpha_2,\cdots, \alpha_N$, and that $\widetilde{\boldsymbol{\Sigma}}_N$ has eigenvalues
	$\widetilde \alpha_1>\widetilde \alpha_2>\cdots>\widetilde \alpha_K$
	with respective   {multiplicities $n_1,\cdots,n_K$, satisfies $1\leq K\leq N$ and $n_1+\cdots n_K=N$}. It is worth noting that the two sets of eigenvalues are connected through a complex mapping that depends on the eigen-structure of ${\boldsymbol{\Sigma}}_N$ and the missing probabilities.   {For example,
		if $\bSigma_N = \diag(\alpha_1,\cdots,\alpha_N)$, then $\widetilde{\bSigma}_N= \diag(\theta_1 \alpha_1,  \cdots, \theta_N \alpha_N).$ It can be noted that in this special case, these two sets of eigenvalues are linked by the missing probabilities. In the general case, the relationship between these two sets of eigenvalues is influenced by the eigenvectors, making the relationship more complex. }
	
	Let $\mathbf{y}_1,\mathbf{y}_2,\dots,\mathbf{y}_n$ be $n$ independent samples drawn from the population $\mathbf{y}$. And we denote the data matrix as  $$\mathbf{Y}=(\mathbf{y}_1,\dots,\mathbf{y}_n)=\left(\mathbf{Y}_N^T ,\mathbf{Y}_S^T\right)^T=\mathbcal{D} \circ \mathbf{X}=\mathbcal{D} \circ \left(\bSigma^{1/2}\mathbcal{X}\right),$$ where  $\mathbcal{D}=({d}_{ij})_{p\times n}$ { consists of independent}   {${d}_{ij}\sim {\rm Ber}(\theta_i)$},  $\mathbf{X}=(\mathbf{x}_1,\dots,\mathbf{x}_n)$, and $\mathbcal{X}=(\mathbcal{x}_1,\dots,\mathbcal{x}_n)$.   {It is
		worth mentioning that $\bY$ is a $p\times n$ matrix,  $\mathbf{Y}_N=(\by_{1,N},\cdots,\by_{n,N})$ is of dimension $N\times n$ since $\by_{i,N}\in\mathbb{R}^N$ and similarly $\mathbf{Y}_S=(\by_{1,S},\cdots,\by_{n,S})$ is of dimension $S\times n$ since $\by_{i,S}\in\mathbb{R}^S$. }  Therefore, the sample covariance matrix is given by 	\begin{equation*}
		\mathbf{S}_n =
		\frac{1}{n}\mathbf{Y}\mathbf{Y}^T
		= 
		\begin{pmatrix}
			\frac{1}{n}\mathbf{Y}_N\mathbf{Y}_N^T & \frac{1}{n}\mathbf{Y}_N\mathbf{Y}_S^T \\
			\frac{1}{n}\mathbf{Y}_S\mathbf{Y}_N^T & \frac{1}{n}\mathbf{Y}_S\mathbf{Y}_S^T
		\end{pmatrix}
		\triangleq
		\begin{pmatrix}
			\mathbf{S}_{11} & \mathbf{S}_{12} \\
			\mathbf{S}_{21} & \mathbf{S}_{22}
		\end{pmatrix}.
	\end{equation*}
	
	Next, we outline the model assumptions.
	\begin{description}
		\item[\textbf{Assumption a: }] As $n\rightarrow\infty$,  $c_n=p/n\to c\in (0, \infty) $. 
		\item[\textbf{Assumption b: }] The entries   {$\mathcal{x}_{ij}$ of $\mathbcal{X}$ satisfy $\mathbb{E}(\mathcal{x}_{ij}^4) <\infty$}. 
		\item[\textbf{Assumption c: }]The empirical spectral distribution (ESD) sequence $H_n$ of $\widetilde{\boldsymbol{\Sigma}}_S$ weakly converges to a probability distribution $H$ as $n\to \infty$.   {Let $\|\boldsymbol{\Sigma}\|$ be bounded by a constant.}
		\item[\textbf{Assumption d: }] The   {distant spiked }eigenvalues $\widetilde \alpha_k,1\le k\le K$
		lie outside $\varGamma_{H}$, the support of $H$, and satisfy $\psi'(\widetilde \alpha_k)>0$ for  all $k$ where
		\begin{equation}\label{cal9}
			\psi(\widetilde \alpha)\triangleq\widetilde \alpha\left[1+c\displaystyle \int \dfrac{t}{\widetilde \alpha - t} \, H(dt) \right]\quad {\rm for}\ \widetilde \alpha \notin \varGamma_{H},\ \widetilde \alpha \neq0.
		\end{equation}
	\end{description} 
	By   {Theorem 2.1 in} \cite{li2024spectral}, the ESD of $\mathbf{S}_{22}$ converges   {almost surely } to $F_{c,H}$ whose Stieltjes transform $m(z)$ satisfies
	\begin{align}\label{m}
		m(z)=\int\frac1{t(1-c-czm(z))-z}H(dt)\ {\rm for}\ z\in\mathbb{C}^+.
	\end{align} 
	Then
	$\underline{F}_{c,H}=(1-c)I_{[0,\infty)}+cF_{c,H}$ is the limiting spectral distribution (LSD) of  $\underline{\mathbf{S}}_{22}=\frac{1}{n}\mathbf{Y}_S^T\mathbf{Y}_S$ and the Stieltjes transform of $\underline{F}_{c,H}$ is denoted by $\underline m(z)$. 	For $\widetilde{\boldsymbol{\Sigma}}_N$, there exists an orthogonal matrix $\mathbf{\widetilde V}=\left(\bv_1,\cdots,   {\bv_N}\right)$ such that \begin{align}\label{se}
		\widetilde{\boldsymbol{\Sigma}}_N=\mathbf{\widetilde V}\rm{diag}\left\{\widetilde \alpha_1{\bf I}_{n_1} ,\dots,\widetilde \alpha_{K}{\bf I}_{n_K}\right\}\mathbf{\widetilde V}^T.
	\end{align}
	
	Let
	$\lambda_{n,1} \ge \lambda_{n,2} \ge \cdots \ge \lambda_{n,p}$
	be the eigenvalues of the sample covariance matrix $\mathbf{S}_n$. The eigenvalues $\lambda_{n,j}, 1 \le j \le N$, are referred to as the sample spiked eigenvalues, with the corresponding eigenvectors termed the sample spiked eigenvectors.
	
	We begin by presenting the convergence of the sample spiked eigenvalues.
	\begin{proposition}\label{pro1}
		Under assumptions (a)-(d), let $\mathcal{k}$ denote the set of indices corresponding to eigenvalues $\widetilde\alpha_k$ with multiplicities $n_k$. Then, we have: 
		\begin{align*}
			  {\lambda_{n,j}-\psi(\widetilde \alpha_k)\xrightarrow{a.s.}0, \quad k=1,\cdots,K,\ j\in\mathcal{k}.}
		\end{align*}
		
		  {Under the assumptions (a)-(c), and for $\tilde{\alpha}_k$ lying outside $\Gamma_H$ with \(\psi'(\widetilde{\alpha}_k)\leq 0\), we have $\lambda_{n,j}\xrightarrow{\text{a.s.}}\sup \Gamma_{F_{c,H}}$ or $\inf \Gamma_{F_{c,H}},$ where $\Gamma_{F_{c,H}}$ is the support of $F_{c,H}$. }
	\end{proposition}
	
	\begin{remark}
		  {In the no-missing  (\(\theta_i = 1\) for all i) and identity-bulk case (\(\boldsymbol{\Sigma}_S = \mathbf{I}_S\)), substituting these conditions into our definition of  $\psi(\widetilde{\alpha})$ recovers the classical spike-to-sample eigenvalue mapping} \[  {\psi(\alpha)=\alpha\left(1+\frac{c}{\alpha-1}\right),}\]
		  {and in this case, the spike phase transition threshold reduces from \(\psi'(\widetilde \alpha)>0\) to \( \alpha > 1+\sqrt{c}\ {\rm or}\  \alpha < 1-\sqrt{c}\), which aligns with the BBP phase transition in \cite{baik2005phase}.}
\end{remark}

\begin{remark}
	The almost surely limit of the sample spiked eigenvalues follows a pattern similar to that of the complete data case. However, the missing data mechanism significantly impacts the results by altering the limiting spectral distribution of the sample covariance matrix for the non-spiked part. { We have some examples to illustrate this below Theorem \ref{thm3}. }
\end{remark}

{
	\begin{remark}
		We consider some simple population models to see how missing schemes influence the spikes. 
		
		\textit{Examples 1.}
		When $\theta_i = \theta$ for $i=1,\cdots, p$, then $\tilde{\bSigma} = \theta^2 \bSigma + (\theta-\theta^2)\bI_p \circ \bSigma = \theta^2 \bSigma + (\theta-\theta^2)\diag(\bSigma)$.  If $\text{diag}(\bSigma) =   {\bI_p}$, $\tilde{\bSigma}=\theta^2 \bSigma + (\theta-\theta^2)\bI_p$. We can see that even a uniform missing scheme for all elements does not scale the eigenvalues proportionally but can have the effect of transferring a spike into a non-spike.
		
		\textit{Examples 2.}
		If $\bSigma = \diag(\alpha_1,\cdots,\alpha_N, 1,\cdots,1)$, then $\tilde{\bSigma}= \diag(\theta_1 \alpha_1,  \cdots, \theta_N \alpha_N, \theta_{N+1}, \cdots, \theta_p).$ For this case, the missing scheme scales the variance of each variable as expected.
	\end{remark}
	
}

Then, we present the  CLT for spiked eigenvalues.
\begin{thm}\label{thm1}
	Under the assumptions (a)-(d),    {the distribution of} the $n_k$-dimensional vector 
	$
	\sqrt{n}\{\lambda_{n,j} - \psi(\widetilde \alpha_k) , j\in \mathcal{k} \}
	$
	converges weakly to  the joint distribution of the eigenvalues of the random matrix 
	\[
	\dfrac{1}{1+cf_2\left(\psi(\widetilde \alpha_k) \right)\widetilde \alpha_k}\left[\mathbf{\widetilde V}^T\mathbcal{C}\left(\psi(\widetilde \alpha_k) \right)\mathbf{\widetilde V} \right]_{kk}, 
	\]
	where $\mathbcal{C}\left(\psi(\widetilde \alpha_k) \right)=(\zeta_{ij}^k)$ is an $N\times N$ symmetric Gaussian random matrix with  zero-mean and covariance function
	\begin{align*}
		{\rm Cov}(\zeta_{j\ell }^k,\zeta_{st}^k)
		&=\tau_0\left(\psi(\widetilde \alpha_k) \right)\left\{\mathcal{E}_{j\ell,st}-\left(\widetilde{\boldsymbol{\Sigma}}_N\right)_{j\ell}\left(\widetilde{\boldsymbol{\Sigma}}_N \right)_{st}\right\} \\
		&\quad+ (\nu_0\left(\psi(\widetilde \alpha_k) \right)-\tau_0\left(\psi(\widetilde \alpha_k) \right))\bigg[\left(\widetilde{\boldsymbol{\Sigma}}_N\right)_{jt}\left(\widetilde{\boldsymbol{\Sigma}}_N \right)_{s\ell}+\left(\widetilde{\boldsymbol{\Sigma}}_N\right)_{js}\left(\widetilde{\boldsymbol{\Sigma}}_N \right)_{t\ell}\bigg], 
	\end{align*} 	
	with $\mathcal{E}_{j\ell,st}=\mathbb{E}\left(y_{j1}y_{s1}y_{\ell1}y_{t1} \right)$ and $\left[\mathbf{\widetilde V}^T\mathbcal{C}\left(\psi(\widetilde \alpha_k) \right)\mathbf{\widetilde V} \right]_{kk}$ is the $k$-th diagonal block using a block decomposition induced by \eqref{se}.  
	Here $s_i=n_1+\cdots+n_{i},i\ge1$, $s_0=0$, $f_2(\lambda)=\displaystyle \int \dfrac{x}{(\lambda - x)^2} \, F_{c,H}(dx)$, and
	\begin{align*}
		\tau_0(\lambda)=\left[\int \dfrac{\lambda}{\lambda - x} \, \underline{F}_{c,H}(dx)\right]^2,\quad \nu_0(\lambda)=c  \int \dfrac{  {\lambda^2}}{(\lambda - x)^2} \, F_{c,H}(dx).
	\end{align*}
\end{thm}

\begin{remark}
	The missing data mechanism has a significant impact on the limiting behavior of the
	sample spiked eigenvalues.   {As formalized in Theorem \ref{thm1}, the missing data  mechanism causes  the mixed moment $\mathcal{E}_{j\ell,st}$ of the observed data in the spiked component involving the mixed moment $\mathbb{E}(d_{j1}d_{s1}d_{\ell1}d_{t1})$ of the missing random variables, and it changes the population spiked eigenvector matrix to $\widetilde{\mathbf{V}}$.} 
	  {
		If $\bSigma = \diag(\alpha_1, 1,\cdots,1)$, then $\widetilde{\bSigma}= \diag(\theta_1 \alpha_1, \theta_{2}, \cdots, \theta_p).$ For this case, $\sqrt{n}\left(\lambda_{n,1}-\psi(\widetilde{\alpha}_1)\right)\xrightarrow{\text{d}}N(0,\sigma_1^2)$, where $\sigma_1^2=\tau_0\left(\psi(\widetilde \alpha_1) \right)\theta_1\mathbb{E}(x_{11}^4)+ (2\nu_0\left(\psi(\widetilde \alpha_1) \right)-3\tau_0\left(\psi(\widetilde \alpha_1) \right))\theta_1^2\alpha_1^2$. This case clearly illustrates how the missing probability affects the limiting variance. }
\end{remark}

In the case where $\widetilde{\boldsymbol{\Sigma}}_N$ has no multiple roots, the joint normality of the spiked eigenvalues can be established.
\begin{itemize}
	\item\textbf{Assumption e: } The eigenvalues of $\widetilde{\boldsymbol{\Sigma}}_N$ are all simple, i.e., $K=N$.
\end{itemize}

\begin{thm}\label{thm11}
	Under the assumptions (a)-(e),    {the distribution of} the $N$-dimensional vector 
	$
	\sqrt{n}\{\lambda_{  {n},k} - \psi(\widetilde \alpha_k) ,1\le k\le N\}
	$
	converges weakly to  the joint distribution of the random vector 
	\[
	\left\{\dfrac{1}{1+cf_2\left(\psi(\widetilde \alpha_k) \right)\widetilde \alpha_k}  {\be_k^T}\left[\mathbf{\widetilde V}^T\mathbcal{C}\left(\psi(\widetilde \alpha_k) \right)\mathbf{\widetilde V} \right]  {\be_k}, 1\le k\le N\right\},
	\]
	  {where $\be_k$ is the $k$-th column of the identity matrix $\bI_N$. }
\end{thm}

The following two theorems provide key results regarding the spiked eigenvectors. {Let $\bbu_k$ be the $k$-th eigenvector of $\widetilde{\bSigma}$,   {where $k\leq N$}.}   {Denote $\boldsymbol{u}_k=[\boldsymbol{u}_{k1}^T,\boldsymbol{u}_{k2}^T]^T$. Since \[\widetilde{\boldsymbol{\Sigma}}\boldsymbol{u}_k=	\begin{pmatrix}
		\widetilde{\boldsymbol{\Sigma}}_N & \mathbf{0} \\
		\mathbf{0} & \widetilde{\boldsymbol{\Sigma}}_S
	\end{pmatrix}\begin{pmatrix}
		\boldsymbol{u}_{k1} \\
		\boldsymbol{u}_{k2}
	\end{pmatrix}=\begin{pmatrix}
		\widetilde{\boldsymbol{\Sigma}}_N\boldsymbol{u}_{k1} \\
		\widetilde{\boldsymbol{\Sigma}}_S\boldsymbol{u}_{k2}
	\end{pmatrix}=\begin{pmatrix}
		\widetilde{\alpha}_k\boldsymbol{u}_{k1} \\
		\widetilde{\alpha}_k\boldsymbol{u}_{k2}
	\end{pmatrix} \] 
	and $\widetilde{\alpha}_k$ is an eigenvalue of $\widetilde{\boldsymbol{\Sigma}}_N$ rather than of $\widetilde{\boldsymbol{\Sigma}}_S$, we have that $\widetilde{\boldsymbol{\Sigma}}_S\boldsymbol{u}_{k2}=\widetilde{\alpha}_k\boldsymbol{u}_{k2}$ if and only if $\boldsymbol{u}_{k2}=\boldsymbol{0}_S$.}   {Thus, $\boldsymbol{u}_{k1}=\bv_{k}$ and}
\begin{align}\label{cal10}
	\widetilde{\boldsymbol{\Sigma}} 
	{\bf u}_{k}=\begin{pmatrix}
		\widetilde{\boldsymbol{\Sigma}}_N & \mathbf{0} \\
		\mathbf{0} & \widetilde{\boldsymbol{\Sigma}}_S
	\end{pmatrix}
	\begin{pmatrix}
		\bv_{k} \\
		\mathbf{0}_{S}
	\end{pmatrix}
	=\widetilde \alpha_k    
	\begin{pmatrix}
		\bv_{k} \\
		\mathbf{0}_{S}
	\end{pmatrix}.
\end{align}
Let 
$
\left(\boldsymbol{\beta}_{1},\cdots,\boldsymbol{\beta}_{p}\right){\rm diag}\left(\lambda_{n,1},\cdots,\lambda_{n,p}\right)	\left(\boldsymbol{\beta}_{1},\cdots,\boldsymbol{\beta}_{p}\right)^T
$
denote the spectral decomposition of $\bS_n$, where $	\left(\boldsymbol{\beta}_{1},\cdots,\boldsymbol{\beta}_{p}\right)$ is an orthogonal matrix consisting of the orthonormal eigenvectors of $\bS_n$. This implies 
\begin{equation}\label{cal4}
	\mathbf{S}_n
	\boldsymbol{\beta}_{k} =
	\begin{pmatrix}
		\mathbf{S}_{11} & \mathbf{S}_{12} \\
		\mathbf{S}_{21} & \mathbf{S}_{22}
	\end{pmatrix}
	\begin{pmatrix}
		\boldsymbol{\beta}_{1k} \\
		\boldsymbol{\beta}_{2k}
	\end{pmatrix}
	=\lambda_{n,k}   
	\begin{pmatrix}
		\boldsymbol{\beta}_{1k} \\
		\boldsymbol{\beta}_{2k}
	\end{pmatrix}, 
\end{equation}
where $	\boldsymbol{\beta}_{1k}$ is a $N$-dimensional vector,     {\(\boldsymbol{\beta}_{2k}\) is an $S$-dimensional vector}, and
\begin{equation}\label{cal8}
	\boldsymbol{\beta}_{1k}^T\boldsymbol{\beta}_{1k}+\boldsymbol{\beta}_{2k}^T\boldsymbol{\beta}_{2k}=1. 
\end{equation}
Then, we provide the limit of the squared inner product between the $k$-th spiked eigenvector of the sample covariance matrix and that of $\widetilde{\boldsymbol{\Sigma}}_N$. 
\begin{thm}\label{thm2}
	Under the assumptions (a)-(e), we have
	\begin{equation*}
		\langle \boldsymbol{\beta}_{k},\mathbf{u}_k\rangle^2=\langle \boldsymbol{\beta}_{1k},\mathbf{\widetilde v}_k\rangle^2 \xrightarrow{\text{a.s.}}
		\dfrac{1}{1+c{f_2}(\psi(\widetilde \alpha_k))\widetilde \alpha_{k}}. 
	\end{equation*}
	  {Under the assumptions (a)-(c) and (e), and for $\tilde{\alpha}_k \notin \Gamma_{H}$ with \(\psi'(\widetilde{\alpha}_k)\leq 0\), we have $\langle \boldsymbol{\beta}_{k},\mathbf{u}_k\rangle^2 \xrightarrow{\text{a.s.}}
		0$.}
\end{thm}

Denote $\mathscr{D}_{k}=\dfrac{\widetilde \alpha_k}{\psi(\widetilde \alpha_k)}\sum\limits_{s\neq k}^{}\dfrac{1}{\widetilde \alpha_{s}-\widetilde \alpha_{k}}\bv_s\bv_s^T .$ Next, we will present the CLT for spiked eigenvectors.
\begin{thm}\label{thm3}
	Recall $\boldsymbol{a}_{k}=\boldsymbol{\beta}_{1k}/\|\boldsymbol{\beta}_{1k}\|$. Under the assumptions (a)-(e), we have
	\begin{align*}
		\sqrt{n}(\boldsymbol{a}_k-\bv_k)& \xrightarrow{\text{d}}N(\mathbf{0},{\boldsymbol{\Gamma}}_k), \\
		\sqrt{n}(\mathbf{\widetilde V}^T\boldsymbol{a}_k-\boldsymbol{e}_k) &\xrightarrow{\text{d}}N(\mathbf{0},\bV^T{\boldsymbol{\Gamma}}_k\bV),
	\end{align*} 
	  {where} $\mathscr{C}^k=(\mathscr{C}^k_{ij}),$ $\mathscr{C}^k_{ij}=\sum_{s,t}^{}\widetilde{v}_{ks}\widetilde{v}_{kt}\mathbb{E}(y_{i1}y_{s1}y_{j1}y_{t1})$, and 
	\begin{align*}
		\boldsymbol{\Gamma}_k
		=&\tau_0\left(\psi(\widetilde \alpha_k) \right)\mathscr{D}_{k}\mathscr{C}^k\mathscr{D}_{k}+\dfrac{\Big(\nu_0\left(\psi(\widetilde \alpha_k) \right)-\tau_0\left(\psi(\widetilde \alpha_k) \right)\Big)}{\psi^2(\widetilde \alpha_k)}\sum\limits_{f\neq k}\dfrac{\widetilde \alpha_k^3\widetilde \alpha_{f}}{(\widetilde \alpha_{f}-\widetilde \alpha_{k})^2}\bv_{f}\bv_{f}^T. 
	\end{align*}
\end{thm}

\section{  {Simulation study}}\label{se3}
  {We conducted simulations to validate the approximation accuracy of the proposed theorem. Set
	$\boldsymbol{\Sigma}=\mathbf{D}\text{diag}\{15,10,2,\dots,2,1,\dots,1\}\mathbf{D}^T$ where the number of $2$ is $248$, and
	\begin{equation*}
	\mathbf{D}=		\begin{pmatrix}
		\mathbf{D}_1 & \mathbf{0} \\
		\mathbf{0} & \mathbf{D}_2
	\end{pmatrix}.
	\end{equation*}
	with $\mathbf{D}_1$ and $\mathbf{D}_2$ both being orthogonal matrices, where $\mathbf{D}_1$ is 2$\times$2. We considered six simulation scenarios with a fixed dimension $p=500$, varying only the sample size $n$ and the underlying distribution, listed in sequence as follows:
	(1): $n = 500$, normal distribution,
	(2): $n = 500$, gamma distribution,
	(3): $n = 500$, student's $t$ distribution,
	(4): $n = 1000$, normal distribution,
	(5): $n = 1000$, gamma distribution, 
	(6): $n = 1000$, student's $t$ distribution},   {where the specific parameter information of the Gamma distribution and the $t$ distribution used here are: $\Gamma(4,0.5)-2$ and  $t_{10}/\sqrt{5/4}$.
}

\begin{figure}[htbp]
	\centering
	\begin{minipage}{0.32\textwidth}
		\includegraphics[width=\textwidth]{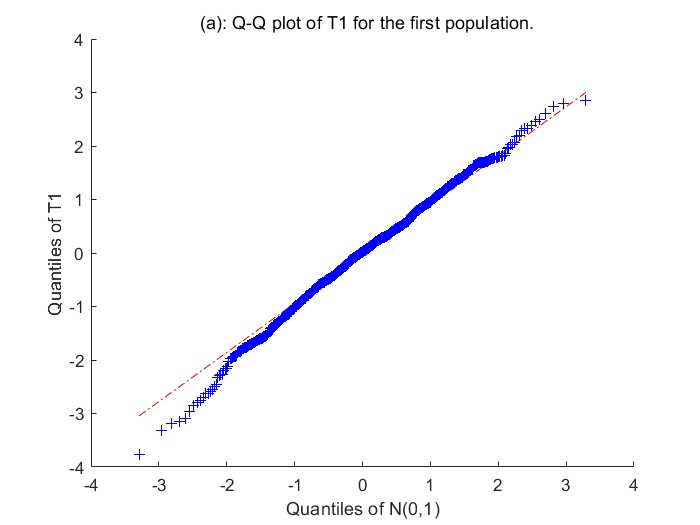}
	\end{minipage}
	\begin{minipage}{0.32\textwidth}
		\includegraphics[width=\textwidth]{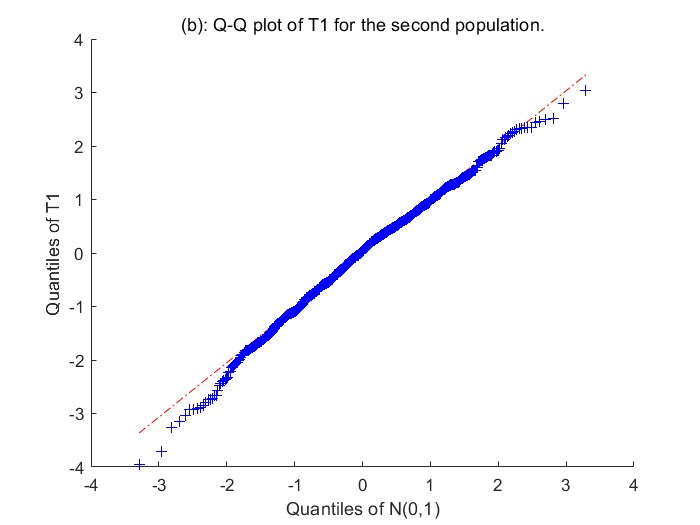}
	\end{minipage}
	\begin{minipage}{0.32\textwidth}
		\includegraphics[width=\textwidth]{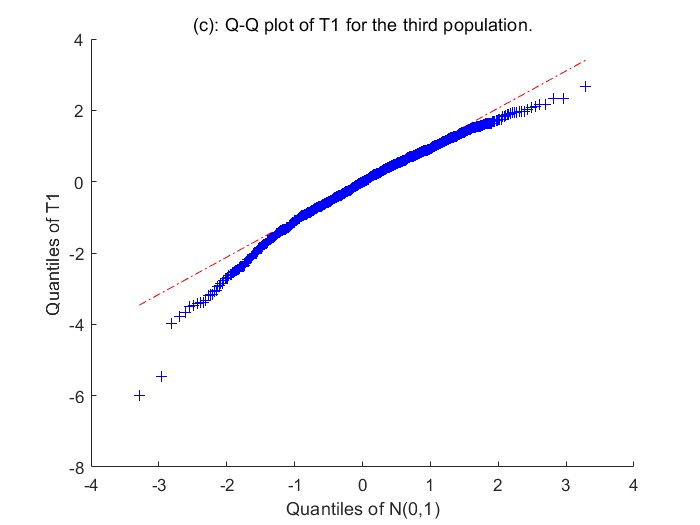}
	\end{minipage}
	\\
	\begin{minipage}{0.32\textwidth}
		\includegraphics[width=\textwidth]{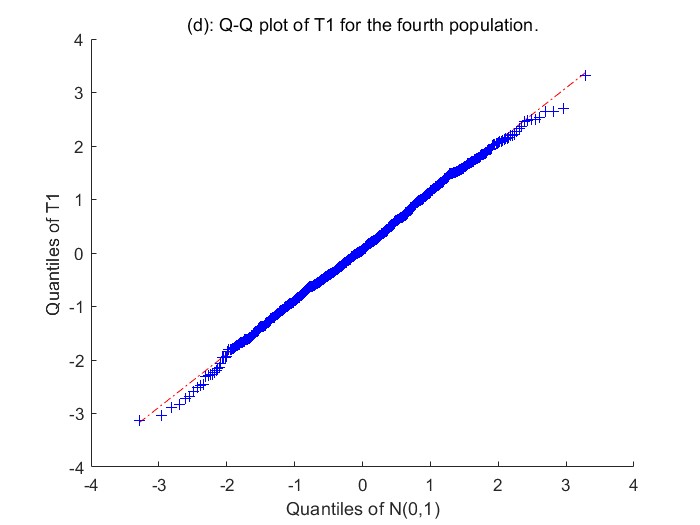}
	\end{minipage}
	\begin{minipage}{0.32\textwidth}
		\includegraphics[width=\textwidth]{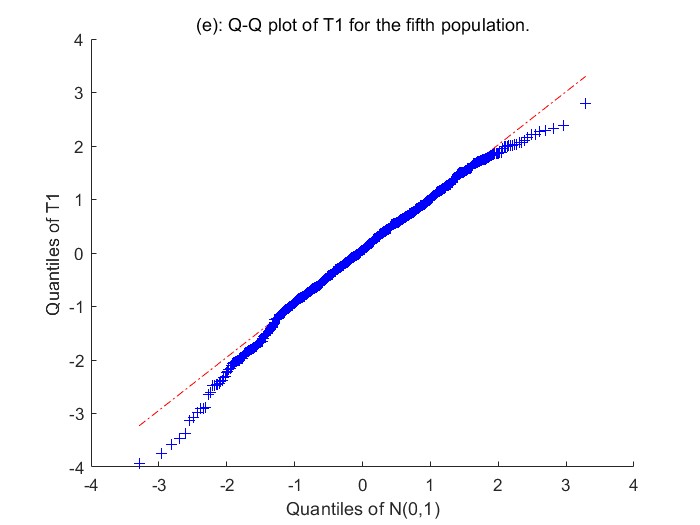}
	\end{minipage}
	\begin{minipage}{0.32\textwidth}
		\includegraphics[width=\textwidth]{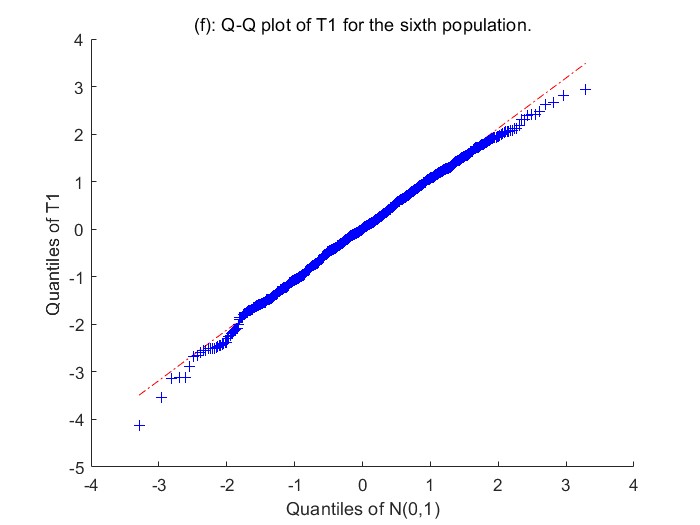}
	\end{minipage}
	\caption{Q-Q plot of the empirical distribution of $T_1$ vs the standard normal distribution. }
	\label{fig:1}
\end{figure}

\begin{figure}[htbp]
	\centering
	\begin{minipage}{0.32\textwidth}
		\includegraphics[width=\textwidth]{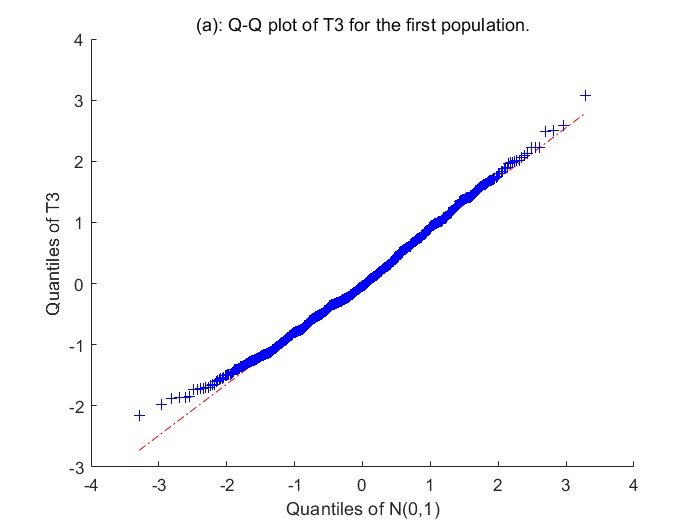}
	\end{minipage}
	\begin{minipage}{0.32\textwidth}
		\includegraphics[width=\textwidth]{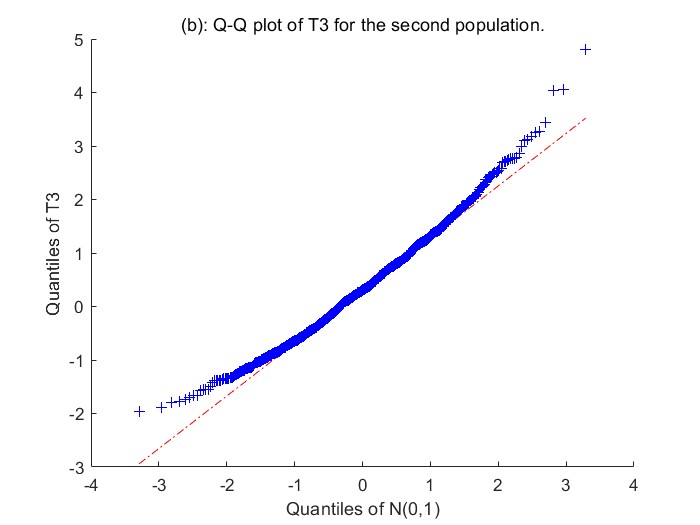}
	\end{minipage}
	\begin{minipage}{0.32\textwidth}
		\includegraphics[width=\textwidth]{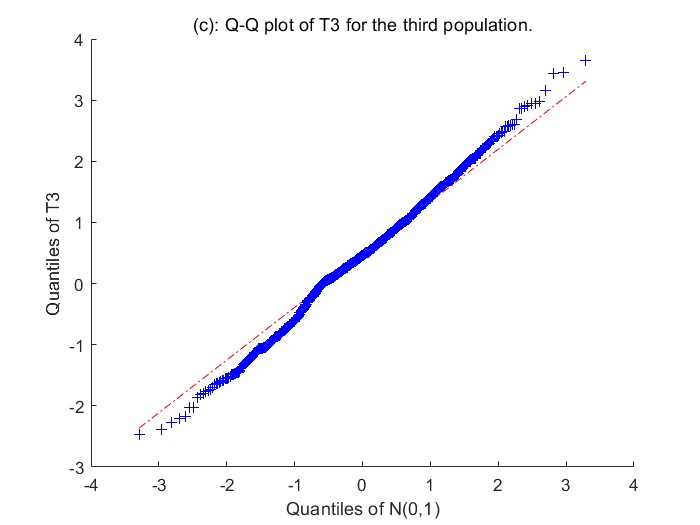}
	\end{minipage}
	\\
	\begin{minipage}{0.32\textwidth}
		\includegraphics[width=\textwidth]{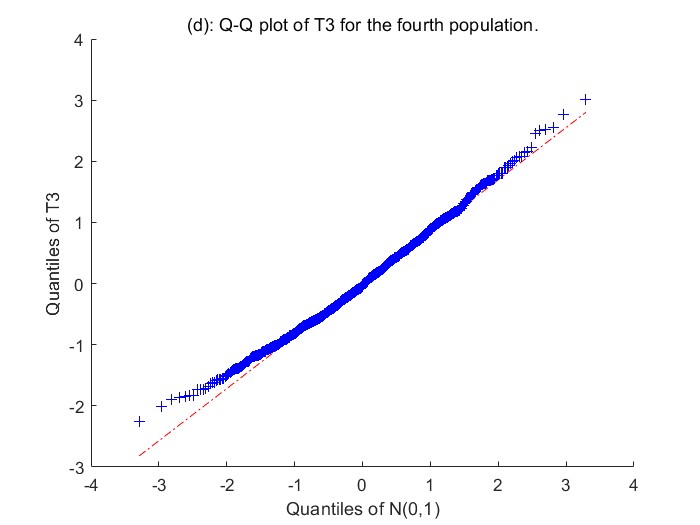}
	\end{minipage}
	\begin{minipage}{0.32\textwidth}
		\includegraphics[width=\textwidth]{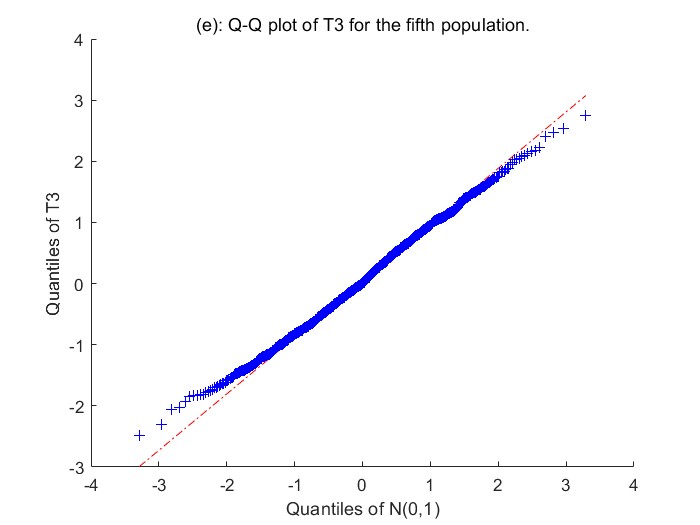}
	\end{minipage}
	\begin{minipage}{0.32\textwidth}
		\includegraphics[width=\textwidth]{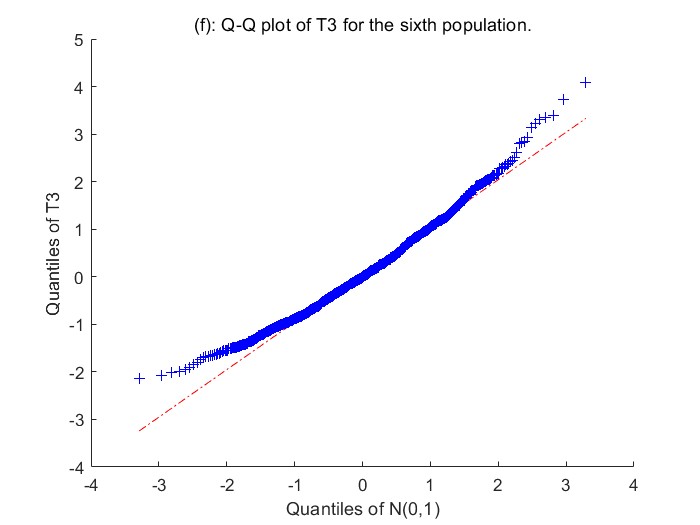}
	\end{minipage}
	\caption{Q-Q plot of the empirical distribution of $T_3$ vs the standard normal distribution. }
	\label{fig:2}
\end{figure}

Among them, the assumed probabilities of missing for each dimension of the population are as follows:
  {$1-\theta_1,\cdots,1-\theta_p\sim{\rm Uniform}(0,0.5)$}. We consider the following statistics:
  {$T_1= \sqrt{n}\{\lambda_{n,1} - \psi(\widetilde \alpha_1) \}/\sqrt{{\rm Cov}(\zeta_{11 }^1,\zeta_{11}^1)},$
	$T_2= \langle \boldsymbol{\beta}_1,\mathbf{u}_1\rangle^2, $ and 
	$ T_3= \sqrt{n}\mathbf{e}_2^T(\mathbf{\widetilde V}^T\boldsymbol{a}_1-\boldsymbol{e}_1)/\sqrt{\be_2^T[\widetilde{\mathbf{V}}^T\boldsymbol{\Gamma}_1\widetilde{\mathbf{V}}]\be_2}. 
	$}
All the above parameters are calculated using Theorem \ref{thm1}, Theorem \ref{thm2}, and Theorem \ref{thm3}. The corresponding simulation results are presented in Figure \ref{fig:1}, Table \ref{tab:1}, and Figure \ref{fig:2}.

\begin{table}[h]
	\centering
	\caption{The simulated values of $T_2$ {  (averaged over 5000 replications)} vs. the theoretical values of $T_2$ }
	\begin{tabular}{cccccccc}
		\toprule
		\midrule
		Value type \textbackslash Population Type & & (1) & (2) & (3) & (4) &(5)&(6)\\
		\midrule
		simulated value    &  &0.8948    & 0.8529 &0.8922 &0.9392  &0.9259&0.9318 \\
		theoretical value   & & 0.9045  & 0.8825&0.8999 & 0.9454 &0.9439&0.9416\\
		\midrule
		\bottomrule
	\end{tabular}
	\label{tab:1}
\end{table}

\section{Application}\label{sec3}
\subsection{An application for testing the independence of two sets of variables}\label{appsec1}
The independence test between two sets of variables is a fundamental aspect of statistical inference, widely applied in various scientific research domains. It is a crucial preliminary step required before conducting numerous scientific analyses. In most cases, the two sets of variables subjected to an independence test represent distinct categories of information. Additionally, in many predictive tasks, these two sets of variables often differ significantly in the amount of information they carry, reflecting their diverse roles in representing different types of information.

If the independence test reveals that the two sets of variables are independent, researchers can design separate neural network branches to process each set of variables individually. The outputs from these branches can then be integrated, allowing the model to leverage the strengths of specialized processing for each feature set. This approach not only enhances the training efficiency of the model but also improves its accuracy.

Therefore, based on Theorem \ref{thm11} and Theorem \ref{thm2}, we present here a method for testing the independence between two sets of variables under the spiked population hypothesis based on high-dimensional missing data. All our assumptions about the data use the data model employed in the section \ref{sec2}. Denote the $p$ variables of the population as $(x_1,\dots,x_p)^T$, and divide the first $N$ variables and the last $S$ variables into two groups of variables $\mathbf{x}_N=(x_1,\dots,x_N)^T$ and $\mathbf{x}_S=(x_{N+1},\dots,x_p)^T$. Therefore, the null hypothesis and alternative hypothesis we consider are as follows:
\[H_0:\ \mathbf{x}_N \ \text{is independent of}\  \mathbf{x}_S \ \ {\rm vs}.\ \ H_1:\ \mathbf{x}_N \ \text{is not independent of}\  \mathbf{x}_S. \]
  {Under $H_0$, the joint distribution of the spikes given in Theorem \ref{thm11} is asymptotically normal. We present the statistic $\mathbcal{T}_N$ that follows a chi-squared distribution with $N$ degrees of freedom in the theorem below. } 
  {\begin{thm}\label{th4.1}
		Under $H_0$, we have $\mathbcal{T}_N:=\mathbcal{H}_N^T\mathbcal{H}_N\sim\chi^2(N)$,
		where 
		\[\mathbcal{H}_N=\boldsymbol{\Sigma}_{1N}^{-1/2}(\sqrt{n}(\lambda_{n,1}-\psi(\widetilde{\alpha}_1)),\dots,\sqrt{n}(\lambda_{n,N}-\psi(\widetilde{\alpha}_N)))^T,\] 
		\begin{align*}
			(\boldsymbol{\Sigma}_{1N})_{ij}=(\mathbf{\widetilde v}_j^T\otimes\mathbf{\widetilde v}_i^T)\frac{1}{1+cf_2(\psi(\widetilde{\alpha}_j))\widetilde{\alpha}_j}\mathbb{E}\left[\mathbcal{C}(\psi(\widetilde{\alpha}_j))\otimes\mathbcal{C}(\psi(\widetilde{\alpha}_i))\right]\frac{1}{1+cf_2(\psi(\widetilde{\alpha}_i))\widetilde{\alpha}_i}  {(\mathbf{\widetilde v}_j\otimes\mathbf{\widetilde v}_i)},
		\end{align*}
		\begin{align*}
			{\rm Cov}\left(\zeta_{k\ell }^j,\zeta_{st}^i\right)&=\widetilde{\tau}_{ji}\left\{\mathcal{E}_{k\ell,st}-\left(\widetilde{\boldsymbol{\Sigma}}_N\right)_{k\ell}\left(\widetilde{\boldsymbol{\Sigma}}_N\right)_{st}\right\} 
			+ (\widetilde{\nu}_{ji}-\widetilde{\tau}_{ji})\left[\left(\widetilde{\boldsymbol{\Sigma}}_N\right)_{kt}\left(\widetilde{\boldsymbol{\Sigma}}_N\right)_{s\ell}+ \left(\widetilde{\boldsymbol{\Sigma}}_N\right)_{ks}\left(\widetilde{\boldsymbol{\Sigma}}_N\right)_{\ell t}  \right] 
			, 
		\end{align*}
		and
		\[\widetilde{\tau}_{ji}=\left[\int \dfrac{\psi(\widetilde{\alpha}_j)}{\psi(\widetilde{\alpha}_j) - x} \, \underline{F}_{c,H}(dx)\right]\cdot\left[\int \dfrac{\psi(\widetilde{\alpha}_i)}{\psi(\widetilde{\alpha}_i) - x} \, \underline{F}_{c,H}(dx)\right], \
		\widetilde{\nu}_{ji}=c  \int \dfrac{x^2F_{c,H}(dx)}{(\psi(\widetilde{\alpha}_j) - x)(\psi(\widetilde{\alpha}_i) - x)} \, . \] 
\end{thm}}
  {\begin{Proof}
		Under $H_0$, according to Theorem \ref{thm11}, $(\sqrt{n}(\lambda_{n,1}-\psi(\widetilde{\alpha}_1)),\dots,\sqrt{n}(\lambda_{n,N}-\psi(\widetilde{\alpha}_N)))^T$ converges weakly to the random vector 
		\[
		\left\{\dfrac{1}{1+cf_2\left(\psi(\widetilde \alpha_k) \right)\widetilde \alpha_k}  {\be_k^T}\left[\mathbf{\widetilde V}^T\mathbcal{C}\left(\psi(\widetilde \alpha_k) \right)\mathbf{\widetilde V} \right]  {\be_k}, 1\le k\le N\right\},
		\]
		where $\mathbcal{C}\left(\psi(\widetilde \alpha_k) \right)=(\zeta_{ij}^k)$ is an $N\times N$ symmetric Gaussian random matrix with  zero-mean and covariance function is known. Simple calculations yield $$\boldsymbol{\Sigma}_{1N}^{-1/2}(\sqrt{n}(\lambda_{n,1}-\psi(\widetilde{\alpha}_1)),\dots,\sqrt{n}(\lambda_{n,N}-\psi(\widetilde{\alpha}_N)))^T\xrightarrow{{\rm d}}N(\boldsymbol{0},\bI_N),$$
		where \begin{align*}
			(\boldsymbol{\Sigma}_{1N})_{ij}=(\mathbf{\widetilde v}_j^T\otimes\mathbf{\widetilde v}_i^T)\frac{1}{1+cf_2(\psi(\widetilde{\alpha}_j))\widetilde{\alpha}_j}\mathbb{E}\left[\mathbcal{C}(\psi(\widetilde{\alpha}_j))\otimes\mathbcal{C}(\psi(\widetilde{\alpha}_i))\right]\frac{1}{1+cf_2(\psi(\widetilde{\alpha}_i))\widetilde{\alpha}_i}(\mathbf{\widetilde v}_j\otimes\mathbf{\widetilde v}_i).
		\end{align*}
		Thus, we can obtain $\mathbcal{T}_{N}\xrightarrow{{\rm d}}\chi^2(N)$, where
		\begin{align*}
		\mathbcal{T}_N={}&(\sqrt{n}(\lambda_{n,1}-\psi(\widetilde{\alpha}_1)),\dots,\sqrt{n}(\lambda_{n,N}-\psi(\widetilde{\alpha}_N)))\boldsymbol{\Sigma}_{1N}^{-1}\\
		&\quad\times(\sqrt{n}(\lambda_{n,1}-\psi(\widetilde{\alpha}_1)),\dots,\sqrt{n}(\lambda_{n,N}-\psi(\widetilde{\alpha}_N)))^T.\qedhere
		\end{align*}
\end{Proof}}

Therefore, for a given significance level $\alpha$, the interval in which we reject $H_0$ based on the statistic   {$\mathbcal{T}_N$ is as follows: 
	$\mathbcal{W}=\left\{\mathbcal{T}_N\geq \chi_{1-\alpha}^2(N)\right\}, $
	where $\chi^2_{1-\alpha}(N)$ is the $1-\alpha$ quantile of the chi-squared distribution with $N$ degrees of freedom, that is, if $\mathbcal{T}_N\sim \chi^2(N)$, then $P(\mathbcal{W})=\alpha$.}  

  {In particular, when the hypothesis \(H_1\) corresponds to the scenario that exactly $L$ components of \(\bx_S\) are dependent on \(\bx_N\) (where $L$ is a fixed positive integer), we present the following theorem to conduct a power analysis of this statistic.}
  {\begin{thm}\label{th4.2}(Power analysis of $\mathbcal{T}_N$)
		Under $H_1$, as $n\rightarrow\infty$, we have $P\left(\mathbcal{T}_N\geq \chi^2_{1-\alpha}(N)\right)\rightarrow1.$
	\end{thm}
	\begin{Proof}
		Under $H_1$, i.e., when the  population covariance matrix satisfies
		\begin{equation*}
		\boldsymbol{\Sigma} =
		\begin{pmatrix}
			\boldsymbol{\Sigma}_{N+L} & \mathbf{0} \\
			\mathbf{0} & \boldsymbol{\Sigma}_{S-L}
		\end{pmatrix} 
		\end{equation*}
		and
		\begin{equation*}
		\widetilde{\boldsymbol{\Sigma}} =
		\begin{pmatrix}
			\widetilde{\boldsymbol{\Sigma}}_{N+L} & \mathbf{0} \\
			\mathbf{0} & \widetilde{\boldsymbol{\Sigma}}_{S-L}
		\end{pmatrix}
		=\mathbcal{P}\boldsymbol{\Sigma}\mathbcal{P}+(\mathbcal{P}-\mathbcal{P}^2)\circ\boldsymbol{\Sigma}, \ {\rm{where}}	\
		\mathbcal{P} =
		\rm{diag}\left\{\theta_1,\dots,\theta_p\right\}
		=
		\begin{pmatrix}
			\mathbcal{P}_{N+L} & \mathbf{0} \\
			\mathbf{0} & \mathbcal{P}_{S-L}
		\end{pmatrix}.
		\end{equation*}
		It can be found by examining the proof of Proposition \ref{pro1} that we can still derive the result $\lambda_{n,k}\xrightarrow{{\rm a.s.}}\psi(\widetilde{\alpha}^{[N+L]}_k)$ in this case,  for $1\leq k\leq N$, where $\widetilde{\alpha}^{[N+L]}_k$ denotes the $k$-th largest eigenvalue of the matrix $\widetilde{\bSigma}_{N+L}$. Since there exists at least one $\widetilde{\alpha}^{[N+L]}_j$, for $1\leq j\leq N$, such that \(\widetilde{\alpha}^{[N+L]}_j > \widetilde{\alpha}_j\) and \(\psi'(x)>0\) holds for distant spikes, we can thus derive that $\psi(\widetilde{\alpha}^{[N+L]}_j)>\psi(\widetilde{\alpha}_j)$. Therefore, we have $\sqrt{n}(\lambda_{n,j}-\psi(\widetilde{\alpha}_j))\rightarrow\infty$, and then $\mathbcal{T}_N\rightarrow\infty$. Thus, under $H_1$, as $n\rightarrow\infty$, we have $P\left(\mathbcal{T}_N\geq \chi^2_{1-\alpha}(N)\right)\rightarrow1.$
\end{Proof}}

However, in practical applications, $\psi(\widetilde{\alpha}_i)$ and $\boldsymbol{\Sigma}_{1N}$ in the statistic $\mathbcal{T}_N$ cannot be directly calculated,  as they involve the unknown population distant spikes $\widetilde{\alpha}_i$ and the first $N$-dimensions $\mathbf{\widetilde v}_i$ of the population distant spiked eigenvectors.   {For $\widetilde{\alpha}_i$, we can estimate it by  $\lambda_{i}(\bS_{11})$, where $\lambda_{i}(\bS_{11})$ denotes the $i$-th largest eigenvalue of the matrix $\bS_{11}$.} For $\mathbf{\widetilde v}_i$, according to Theorem \ref{thm2} of this paper, the first $N$ dimensions of the $i$-th sample spiked eigenvectors converge almost surely to $\mathbf{\widetilde v}_i$ after normalization. Therefore, we can choose this method to obtain an estimate $\hat{\mathbf{\widetilde v}}_i$ of $\mathbf{\widetilde v}_i$.  So we get   { $\hat{\widetilde{\alpha}}_i=\lambda_{i}(\bS_{11})$ and $\hat{\mathbf{\widetilde v}}_i=\boldsymbol{\beta}_{1i}/\|\boldsymbol{\beta}_{1i}\|  \ \text{ for} \ i=1,\dots,N$.   {And, for $\psi(\widetilde{\alpha}_i)$, we can estimate it by $\hat{\psi}(\hat{\widetilde{\alpha}}_i)$. }According to (2.29) and (2.30) in \cite{zhang2022asymptotic}, we have $\hat{\psi}'(\hat{\widetilde{\alpha}}_i)=\frac{1}{\hat{\widetilde{\alpha}}_i^2\hat{\underline{m}}'\{\psi(\widetilde{\alpha}_i)\}}$, $\hat{\underline{m}}'\{\psi(\widetilde{\alpha}_i)\}=\frac{1-c_n}{\lambda^2_i}+\frac{1}{n}\sum_{j=N+1}\frac{1}{(\lambda_j-\lambda_i)^2}$. Naturally, we can further obtain the estimator:  $\hat{f}_2(\hat{\psi}(\hat{\widetilde{\alpha}}_i))=\frac{1}{p-N}\sum_{j=N+1}^{p}\lambda_j(\hat{\psi}(\hat{\widetilde{\alpha}}_i)-\lambda_j)^{-2}$.}

  {Based on Theorem \ref{th4.1}, the results in Section 2.2 of \cite{waternaux1976asymptotic}, and the delta method, we can present the asymptotic distribution of the estimator $\hat{\mathbcal{T}}_N$ for $\mathbcal{T}_N$ in the following theorem. Furthermore, in practical applications, we use half of the data to obtain the spiked eigenvalues $\lambda_{n,1},\cdots,\lambda_{n,N}$, and the other half to estimate the unknown parameters.}

  {\begin{thm}
		Under $H_0$, we have $\hat{\mathbcal{T}}_N:=\hat{\mathbcal{H}}_N^T\hat{\mathbcal{H}}_N\sim\chi^2(N)$,
		where 
		\[\hat{\mathbcal{H}}_N=\left(\hat{\boldsymbol{\Sigma}}_{1N}+\hat{\Omega}_N\right)^{-1/2}(\sqrt{n/2}(\lambda_{n,1}-\hat{\psi}(\hat{\widetilde{\alpha}}_1)),\dots,\sqrt{n/2}(\lambda_{n,N}-\hat{\psi}(\hat{\widetilde{\alpha}}_N)))^T,\] \[(\hat{\Omega}_N)_{ij}=\hat{\psi}'(\hat{\widetilde{\alpha}}_i)\hat{\psi}'(\hat{\widetilde{\alpha}}_j)(\mathbf{\hat{\widetilde v}}_i^T\mathbb{K}_4\mathbf{\hat{\widetilde v}}_i\times\mathbf{\hat{\widetilde v}}_j^T\mathbb{K}_4\mathbf{\hat{\widetilde v}}_j+2\hat{\widetilde{\alpha}}_i\hat{\widetilde{\alpha}}_jI(i=j)),\]
		\[ (\mathbb{K}_4)_{ij}=\mathbb{E}y_{1i}^2y_{1j}^2-\mathbb{E}y_{1i}^2\mathbb{E}y_{1j}^2-2(\mathbb{E}y_{1i}y_{1j})^2,\]
		and
		\[(\hat{\boldsymbol{\Sigma}}_{1N})_{ij}=(\mathbf{\hat{\widetilde v}}_j^T\otimes\mathbf{\hat{\widetilde v}}_i^T)\frac{1}{1+2c_n\hat{f}_2(\hat{\psi}(\hat{\widetilde{\alpha}}_j))\hat{\widetilde{\alpha}}_j}\mathbb{E}\left[\mathbcal{C}(\hat{\psi}(\hat{\widetilde{\alpha}}_j))\otimes\mathbcal{C}(\hat{\psi}(\hat{\widetilde{\alpha}}_i))\right]\frac{1}{1+2c_n\hat{f}_2(\hat{\psi}(\hat{\widetilde{\alpha}}_i))\hat{\widetilde{\alpha}}_i}(\mathbf{\hat{\widetilde v}}_j\otimes\mathbf{\hat{\widetilde v}}_i).\]
\end{thm}}

  {\begin{remark}
		The population moments involved in the statistic $\hat{\mathbcal{T}}_N$ are all estimated using sample moments.
\end{remark}}

  {Furthermore, if $\hat{\mathbcal{T}}_N\sim \chi^2(N)$, then $P(\hat{\mathbcal{W}})=P(\hat{\mathbcal{T}}_N\geq \chi^2_{1-\alpha}(N))=\alpha$. Similar to the proof of Theorem \ref{th4.2}, we can similarly obtain the following theorem.
	\begin{thm}(Power analysis of $\hat{\mathbcal{T}}_N$)
		Under $H_1$, as $n\rightarrow\infty$, we have $P(\hat{\mathbcal{T}}_N\geq \chi^2_{1-\alpha}(N))\rightarrow1$.
\end{thm}}


\subsection{Numerical results for section \ref{appsec1}}
In this section, we propose a simulation study for the case  where $N=2$ to evaluate the performance of the statistic proposed in Section \ref{appsec1} for testing the independence between two sets of variables under the assumption of the spiked population model in the presence of missing data.   {We also compare it with the statistics $T_{W}$ and $T_{LH}$ given in \cite{bodnar2019testing}. }

We considered populations with dimension-to-sample size ratios $c_n\in\left\{  {0.4,0.6,0.8}\right\}$ and sample sizes \(n \in \{100, 250, 500\}\). For populations with different dimensions, we generate independently and identically distributed data with zero expectation and a covariance matrix of $\boldsymbol{\Sigma}_p$ where $\boldsymbol{\Sigma}_p=\mathbf{D}{\rm diag}\left\{  {0.6,0.4,0.02\times\boldsymbol{1}_{p/2-2}^T,0.01\times\boldsymbol{1}_{p/2}^T}\right\}\mathbf{D}^T$ and 
\begin{equation*}
\mathbf{D} =
\begin{pmatrix}
	\mathbf{D}_1 & \mathbf{0} \\
	\mathbf{0} & \mathbf{D}_2
\end{pmatrix}.
\end{equation*}
When simulating the empirical size, we take $\mathbf{D}_1$ as a 2-dimensional orthogonal matrix and $\mathbf{D}_2$ as a $(p-2)$-dimensional orthogonal matrix.   {For the empirical power simulation, we take $N+L=15$, and the covariance between each pair of the $L$ variables in $\bx_S$ and $\bx_N$ is assigned as 0.5.} For each combination of $p$ and $n$, we consider   {three} different distributions: the Gaussian distribution,   the Gamma distribution and   {the $t$ distribution}. 

Afterwards, we set a random vector $\mathbcal{b}=(b_1,\dots,b_p)^T$ where each dimension follows a Bernoulli distribution where   {$b_i\sim {\rm Ber}(\theta_i)$} for $i=1,\dots,p$. And during the simulation, after setting 
  {$\theta_1,\cdots,\theta_p\sim{\rm Uniform}(0.5,1)$}, we take the hadamard product of the matrix obtained by sampling the random vector $\mathbcal{b}$ for $n$ times and the previously generated data matrix to obtain our matrix of missing data. Each simulation experiment is repeated 10000 times based on the matrix of missing data generated in this way to obtain the empirical size and the empirical power.

\begin{table}[htbp]
	\centering 
	\setlength{\tabcolsep}{4pt} 
	\caption{Empirical sizes for tests $\hat{\mathbcal{T}}_N$, $T_W$ and $T_{LH}$ under Gaussian, Gamma, and $t$ distributions. }
	{%
		\begin{tabular}{lcccccccccc}
			\hline
			& & \multicolumn{3}{c}{$Gaussian$} & \multicolumn{3}{c}{$Gamma$} & \multicolumn{3}{c}{$t$} \\
			\cline{3-11} 
			test & $n$ & $c_n=0.4$ & $c_n=0.6$ & $c_n=0.8$ & $c_n=0.4$ & $c_n=0.6$ & $c_n=0.8$ & $c_n=0.4$ & $c_n=0.6$ & $c_n=0.8$  \\
			\hline
			
			& 100 & 4.70 & 4.19 & 4.64 & 3.26 & 3.17 & 3.11 & 5.78 & 4.62 & 5.42 \\
			$\hat{\mathbcal{T}}_N$ 	& 250 & 5.06 & 4.93 & 6.90 & 4.72 & 6.25 & 4.61 & 5.45 & 4.04 & 5.10 \\
			& 500 & 5.54 & 5.86 & 5.58 & 4.30 & 4.82 & 5.84 & 5.74 & 4.87 & 5.97 \\
			
			& 100 & 5.35 & 5.19 & 5.40 & 5.24 & 5.47 & 4.96 & 5.29 & 5.75 & 5.29 \\
			$T_W$ 	& 250 & 5.23 & 5.51 & 5.07 & 5.07 & 4.82 & 5.50 & 5.14 & 5.42 & 5.44\\
			& 500 & 5.30 & 4.86 & 5.10 & 4.98 & 4.75 & 5.58 & 4.93 & 4.92 & 4.55 \\
			
			& 100 & 5.36 & 5.48 & 6.50 & 5.26 & 5.84 & 5.63 & 5.06 & 6.12 & 6.06 \\
			$T_{LH}$	& 250 & 5.17 & 5.55 & 5.41 & 5.28 & 5.16 & 5.27 & 5.17 & 5.61 & 5.95 \\
			& 500 & 5.31 & 4.89 & 5.57 & 5.05 & 4.82 & 5.59 & 4.94 & 4.94 & 4.61 \\
			\hline
		\end{tabular}\label{tab:2}
	}
\end{table}

\begin{table}[htbp]
	\centering 
	\setlength{\tabcolsep}{4pt} 
	\caption{Empirical powers for tests $\hat{\mathbcal{T}}_N$, $T_W$ and $T_{LH}$ under Gaussian, Gamma, and $t$ distributions. }
	\label{tab:reject_prob}
	{%
		\begin{tabular}{lcccccccccc}
			\hline
			& & \multicolumn{3}{c}{$Gaussian$} & \multicolumn{3}{c}{$Gamma$} & \multicolumn{3}{c}{$t$} \\
			\cline{3-11} 
			test & $n$ & $c_n=0.4$ & $c_n=0.6$ & $c_n=0.8$ & $c_n=0.4$ & $c_n=0.6$ & $c_n=0.8$ & $c_n=0.4$ & $c_n=0.6$ & $c_n=0.8$  \\
			\hline
			
			& 100 & 69.29 & 77.41 & 79.50 & 59.89 & 65.08 & 67.40 & 59.40 & 71.71 & 74.27 \\
			$\hat{\mathbcal{T}}_N$ 	& 250 & 90.06 & 97.57 & 95.12 & 85.24 & 89.95 & 89.07 & 89.76 & 90.99 & 92.65 \\
			& 500 & 98.41 & 98.71 & 98.43 & 92.74 & 93.45 & 93.07 & 94.27 & 95.01 & 95.42 \\
			
			& 100 & 12.16 & 7.00 & 2.83 & 19.76 & 7.83 & 3.06 & 11.50 & 8.06 & 2.90 \\
			$T_W$ 	& 250 & 41.30 & 10.39 & 6.61 & 44.67 & 19.10 & 4.48 & 33.89 & 20.57 & 3.42 \\
			& 500 &65.00 & 37.61 & 11.48 & 35.06 & 40.84 & 12.43 & 69.52 & 38.88 & 11.92 \\
			
			& 100 & 19.97 & 11.80 & 5.36 & 31.71 & 13.08 & 5.64 & 18.96 & 12.87 & 5.27 \\
			$T_{LH}$ 	& 250 & 56.64 & 15.15 & 10.51 & 59.83 & 27.14 & 6.93 & 47.92 & 28.59 & 5.47 \\
			& 500 & 78.16 & 48.13 & 15.91 & 45.55 & 51.42 & 17.03 & 82.26 & 49.56 & 16.21 \\
			\hline
		\end{tabular}\label{tab:3}
	}
\end{table}

  {Table \ref{tab:2} presents the empirical sizes of the statistics \(\hat{\mathbcal{T}}_N\), $T_W$, and $T_{LH}$, calculated under three specific data distributions (Gaussian distribution, Gamma distribution, and Student’s t-distribution) at the significance level \(\alpha = 0.05\). Correspondingly, Table \ref{tab:3} reports the empirical powers of these same statistics 
	for the identical three data distributions at this 0.05 significance level. }

It can be seen from the simulation results that the testing method we proposed performs well in the
independence test of two sets of variables with high-dimensional missing data.   {Furthermore, simulations of the empirical size and empirical power of our proposed statistic under different missing probabilities reveal that its performance does not exhibit significant variations across different missing probabilities (these simulation results are provided in the supplementary material). Only when the data are excessively missing will the performance be affected due to the failure to estimate the moments of random variables. And, excessively high missing probabilities may cause
	spikes to convert into non-spikes.} 
\section{Proof of theoretical results}\label{sec4}
\subsection{Proof of Proposition \ref{pro1}, Theorem \ref{thm1}  and Theorem \ref{thm11}}\label{sec11}

We present only the proof of Proposition \ref{pro1} and Theorem \ref{thm1}, as the proof of Theorem \ref{thm11} is similar and therefore omitted.
\subsubsection{  {Proof of Proposition \ref{pro1}}}
For any   {distant} spiked eigenvalue $\lambda$ of $\mathbf{S}_n$, we have 
\begin{align*}
	0&=|\lambda \mathbf{I}_{  {p}}-\mathbf{S}_n|=|\lambda \mathbf{I}_{  {S}}-\bS_{22}|\cdot|\lambda\mathbf{I}_{  {N}}-\bS_{11}-\bS_{12}(\lambda \mathbf{I}_{  {S}}-\bS_{22})^{-1}\bS_{21}|,
\end{align*}
  {where \(|\cdot|\) denotes the determinant}. 

By   {Theorem 2.2 } in \cite{li2024spectral},   {when $n\rightarrow\infty$, with probability 1, the eigenvalues of $\bS_{22}$ lie within $\varGamma_{c,H}$ , the support of $F^{c,H}$, and the probability of the event 
	$\mathscr{H}_n = \{\lambda_{n,j} \notin \varGamma_{c,H}\ \cap\ \text{spectrum of} \ \bS_{22} \subset \varGamma_{c,H}\}$ tends to 1. Conditional on this event, }
we only need to consider the equation:
\begin{align}\label{cal3}
	0&=|\lambda\mathbf{I}_{  {N}}-\bS_{11}-\bS_{12}(\lambda \mathbf{I}_{  {S}}-\bS_{22})^{-1}\bS_{21}|=|\lambda \mathbf{I}_{  {N}}-\mathbf{A}_n(\lambda)|,
\end{align}
where
\begin{align}\label{cal7}
	\mathbf{A}_n(\lambda)=\dfrac{1}{n}\mathbf{Y}_N[\mathbf{I}_{  {n}}+\mathbf{B}_n(\lambda)]\mathbf{Y}_N^T, \qquad
	&\mathbf{B}_n(\lambda)=\dfrac{1}{n}\mathbf{Y}_S^T(\lambda \mathbf{I}_{  {S}}-\dfrac{1}{n}\mathbf{Y}_S\mathbf{Y}_S^T)^{-1}\mathbf{Y}_S.
\end{align}
Noting that $\mathbf{Y}_N$ is independent of $\mathbf{B}_n(\lambda)$, we have the decomposition 
\begin{align}\label{cal1}
	\begin{split}
		\mathbf{A}_n(\lambda)&=\dfrac{1}{n}\left\{\mathbf{Y}_N[\mathbf{I}_{  {n}}+\mathbf{B}_n(\lambda)]\mathbf{Y}_N^T- \widetilde{\boldsymbol{\Sigma}}_N{\rm{tr}} (\mathbf{I}_{  {n}}+\mathbf{B}_n(\lambda))\right\}+\dfrac{1}{n} \widetilde{\boldsymbol{\Sigma}}_N{\rm{tr}}(\mathbf{I}_{  {n}}+\mathbf{B}_n(\lambda))\\
		&=\dfrac{1}{\sqrt{n}}\mathbf{C}_n(\lambda)+[1+  {\frac{S}{n}}f_1(\lambda)]\widetilde{\boldsymbol{\Sigma}}_N+  {\frac{S}{n}}\big(  {\dfrac{1}{S}} {\rm{tr}}\mathbf{B}_n(\lambda)-f_1(\lambda)\big)\widetilde{\boldsymbol{\Sigma}}_N, \end{split}
\end{align}
where $f_1(\lambda)=\displaystyle \int_{}^{} \dfrac{x}{\lambda - x} \, F_{c,H}(dx)$ and
$
\mathbf{C}_n(\lambda)=\dfrac{1}{\sqrt{n}}\left\{\mathbf{Y}_N[\mathbf{I}_{  {n}}+\mathbf{B}_n(\lambda)]\mathbf{Y}_N^T- \widetilde{\boldsymbol{\Sigma}}_N{\rm{tr}}\left(\mathbf{I}_{  {n}}+\mathbf{B}_n(\lambda)\right)\right\}. 
$
Due to the almost sure convergence of the Stieltjes transform of the ESD of $  {\frac{1}{n}}\mathbf{Y}_S\mathbf{Y}_S^T$, which is proved in   {Theorem 2.1 of } \cite{li2024spectral}, we know that 
\begin{align*}
	  {\frac{S}{n}}\bigg(  {\frac{1}{S}} \mathrm{tr}\mathbf{B}_n(\lambda) - f_1(\lambda)\bigg) =   {\frac{S}{n}\bigg(\frac{1}{S} \mathrm{tr}\left[(\lambda \mathbf{I}_{  {S}}-\frac{1}{n}\mathbf{Y}_S\mathbf{Y}_S^T)^{-1}\dfrac{1}{n}\mathbf{Y}_S\mathbf{Y}_S^T\right] - f_1(\lambda)\bigg)}= o_{a.s.}(1).
\end{align*}
This, combined with the moment convergence of $\mathbf{C}_n(\lambda)$ proved below, ensures the almost sure convergence of $\mathbf{A}_n(\lambda)-[1 + c f_1(\lambda)]\widetilde{\bSigma}_N$ to $0$. 
  {Thus, for $j\in\mathcal{k}$, we have the determinant of the following $N\times N$ matrix: 
	$$\bV\begin{pmatrix}
		(\lambda_{n,j}-[1+cf_1(\lambda_{n,j})]\widetilde \alpha_1+o_{{\rm a.s.}}(1))\mathbf{I}_{n_1} &&&\\
		&&\ddots &\\
		&&&(\lambda_{n,j}-[1+cf_1(\lambda_{n,j})]\widetilde \alpha_K+o_{{\rm a.s.}}(1))\mathbf{I}_{n_K}
	\end{pmatrix}\bV^T$$
	equal to zero, which is also to say that $\lambda_{n,j}$ satisfies the equation:
	\begin{align*}
		&\lambda_{n,j}-[1+cf_1(\lambda_{n,j})]\widetilde \alpha_k+o_{{\rm a.s.}}(1)\\
		=&\lambda_{n,j}-\widetilde{\alpha}_k\left(1-c-c\lambda_{n,j}m(\lambda_{n,j})\right)+o_{{\rm a.s.}}(1)=\lambda_{n,j}+\lambda_{n,j}\underline{m}(\lambda_{n,j})\widetilde{\alpha}_k+o_{{\rm a.s.}}(1)=0. 
	\end{align*}
	Combining this with (\ref{m}), we have
	\begin{align*}
		&\lambda_{n,j}-\widetilde{\alpha}_k\left(1-c-c\lambda_{n,j}m(\lambda_{n,j})\right)+o_{{\rm a.s.}}(1)\\
		=&\lambda_{n,j}-\widetilde \alpha_k\left[1+c\displaystyle \int \dfrac{t}{\widetilde \alpha_k - t} \, H(dt) \right]+o_{{\rm a.s.}}(1)\triangleq\lambda_{n,j}-\psi(\widetilde \alpha_k)+o_{{\rm a.s.}}(1)=0.
\end{align*}}

{ Specifically, when \(\psi'(\widetilde{\alpha}_k)\leq 0\), leveraging the exact separation results under missing observations established in Section 2.5 of \cite{li2024spectral}, we can analogously derive $\lambda_{n,j}\xrightarrow{\text{a.s.}}\sup \Gamma_{F_{c,H}}$ or $\inf \Gamma_{F_{c,H}},$ via the method proposed in \cite{bai2012sample}, where $\Gamma_{F_{c,H}}$ denotes the support of $F_{c,H}$.}

The proof of Proposition \ref{pro1} is complete.\hfill$\square$
\subsubsection{  {Proof of Theorem \ref{thm1}}}
  {By Theorem 2.3} in \cite{li2024spectral}, we have 
$
  {\frac{S}{n}}\bigg(  {\frac{1}{S}} \mathrm{tr}\mathbf{B}_n(\lambda) - f_1(\lambda)\bigg) = o_{p}\bigg(\frac{1}{\sqrt{n}}\bigg).
$
Thus, to prove Theorem \ref{thm1}, our goal is to find the limiting distribution of $\mathbf{C}_n(\lambda)$.

Let the $i$-th row of $\mathbf{Y}_N$ be denoted by $\by_{(i)}$. It   {can be found} that each element of  $\mathbf{C}_n(\lambda)$ can be represented as
\begin{align}\label{cal2}
	\dfrac{1}{\sqrt{n}}\left[\by_{(i)}(\mathbf{I}_n+\mathbf{B}_n(\lambda))\by_{(j)}^T- (\widetilde{\boldsymbol{\Sigma}}_N)_{ij}{\rm{tr}}(\mathbf{I}_n+\mathbf{B}_n(\lambda))\right],\ i,j=1,\cdots,N.\end{align}
Considering $\mathbf{C}_n(\lambda)$ is a $N \times N$ symmetric random matrix,  it suffices to establish the limiting joint distribution of $g={N(N+1)}/{2}$ random quadratic variables derived from $\mathbf{C}_n(\lambda)$ as specified in \eqref{cal2}.  By Lemma \ref{Lemma 1: } and Lemma \ref{Lemma 2:}, we can obtain $\text{vech}\left(\mathbf{C}_n(\lambda)\right)\xrightarrow{\text{d}} N_g(\mathbf{0},\mathbf{C}).$
Here, $(\mathbf{C})_{ij}$ is the covariance between the $i$-th and $j$-th random quadratic forms of $\text{vech}(\mathbf{C}_n(\lambda))$. 
Therefore, we see that the asymptotic distribution of $\mathbf{C}_n(\lambda)$ is the same as the distribution of $\boldsymbol{\mathbcal{C}}(\lambda)$, whose elements follow normal distributions satisfying $\mathbb{E}(\boldsymbol{\mathbcal{C}}(\lambda)) = 0$ and ${\rm Cov}\left((\bmC(\lambda))_{j\ell},(\bmC(\lambda))_{st}\right)$ being
\begin{align*}
	&\tau_0(\lambda)\left\{\mathbb{E}\left(d_{j1}d_{s1}d_{\ell1}d_{t1}x_{j1}x_{s1}x_{\ell1}x_{t1} \right)-\left(\widetilde{\boldsymbol{\Sigma}}_N\right)_{j\ell}\left(\widetilde{\boldsymbol{\Sigma}}_N \right)_{st}\right\} \\
	&+ (\nu_0(\lambda)-\tau_0(\lambda))\left[\left(\widetilde{\boldsymbol{\Sigma}}_N\right)_{jt}\left(\widetilde{\boldsymbol{\Sigma}}_N \right)_{s\ell}+ \left(\widetilde{\boldsymbol{\Sigma}}_N\right)_{js}\left(\widetilde{\boldsymbol{\Sigma}}_N \right)_{t\ell}\right].
\end{align*}

Let $\chi_{k,j}=\sqrt{n}(\lambda_{n,j}-\psi(\widetilde \alpha_k)),\  j\in \mathcal{k}$. Similar to what was done in section 6.4 of \cite{bai2008central}, we have the following decomposition
\begin{equation*}
	\lambda_{n,j}\mathbf{I}_{  {N}}-\mathbf{A}_n(\lambda_{n,j})=\psi(\widetilde \alpha_k)\mathbf{I}_{  {N}}+\frac{1}{\sqrt{n}}\chi_{k,j}\mathbf{I}_{  {N}}-\mathbf{A}_n(\psi(\widetilde \alpha_k))-\left[\mathbf{A}_n(\lambda_{n,j})-\mathbf{A}_n(\psi(\widetilde \alpha_k))\right], 
\end{equation*}
where
$
\mathbf{A}_n(\lambda_{n,j})-\mathbf{A}_n(\psi(\widetilde \alpha_k))=-\dfrac{1}{\sqrt{n}}\chi_{k,j}\left[cf_2(\psi(\widetilde \alpha_k))\widetilde{\boldsymbol{\Sigma}}_N+o_p(1)\right].
$
Together with \eqref{cal1}, it follows that
\begin{align*}
	\lambda_{n,j}\mathbf{I}_{  {N}}-\mathbf{A}_n(\lambda_{n,j})&=\psi(\widetilde \alpha_k)\mathbf{I}_{  {N}}+\dfrac{1}{\sqrt{n}}\chi_{k,j}\left[\mathbf{I}_{  {N}}+cf_2(\psi(\widetilde \alpha_k))\widetilde{\boldsymbol{\Sigma}}_N\right]-\dfrac{1}{\sqrt{n}}\mathbf{C}_n(\psi(\widetilde \alpha_k))\\
	&-\left[1+cf_1(\psi(\widetilde \alpha_k))\right]\widetilde{\boldsymbol{\Sigma}}_N+o_p(\frac{1}{\sqrt{n}}).
\end{align*}
Replacing $\lambda$ with $\psi(\widetilde \alpha_k)$, one has
\begin{equation}\label{cal19}
	\mathbf{C}_n(\psi(\widetilde \alpha_k))\xrightarrow{\text{d}}\mathbcal{C}(\psi(\widetilde \alpha_k)).
\end{equation}
Furthermore, it implies by Skorokhod strong representation $\mathbf{C}_n(\psi(\widetilde \alpha_k))\xrightarrow{\text{a.s.}}\mathbcal{C}(\psi(\widetilde \alpha_k))$ on an appropriate probability space.
By \eqref{se}, we deduce
\begin{align*}
	&\mathbf{\widetilde V}^T\left[\lambda_{n,j}\mathbf{I}_{  {N}}-\mathbf{A}_n(\lambda_{n,j})\right]\mathbf{\widetilde V}\\&=\begin{pmatrix}
		(\psi(\widetilde \alpha_k)-[1+cf_1(\psi(\widetilde \alpha_k))]\widetilde \alpha_1)\mathbf{I}_{n_1} &&&\\
		&&\ddots &\\
		&&&(\psi(\widetilde \alpha_k)-[1+cf_1(\psi(\widetilde \alpha_k))]\widetilde \alpha_K)\mathbf{I}_{n_K}
	\end{pmatrix}
	\\&+\dfrac{1}{\sqrt{n}}\chi_{k,j}\begin{pmatrix}
		[1+cf_2(\psi(\widetilde \alpha_k))\widetilde \alpha_1]\mathbf{I}_{n_1} &&&\\
		&&\ddots &\\
		&&&[1+cf_1(\psi(\widetilde \alpha_k))\widetilde \alpha_K]\mathbf{I}_{n_K}
	\end{pmatrix}\\
	&
	-\dfrac{1}{\sqrt{n}}\mathbf{\widetilde V}^T\mathbcal{C}(\psi(\widetilde \alpha_k))\mathbf{\widetilde V}+o_{a.s.}(\dfrac{1}{\sqrt{n}}).
\end{align*}
It can be verified that $
\psi(\widetilde \alpha_k)-[1+cf_1(\psi(\widetilde \alpha_k))]\widetilde \alpha_k=0.
$
Due to \eqref{cal3} and the above equality, we find 
\begin{equation*}
	\left|\dfrac{1}{\sqrt{n}}\chi_{k,j}\left(1+cf_2(\psi(\widetilde \alpha_k))\widetilde \alpha_k\right)\mathbf{I}_{n_k}-\dfrac{1}{\sqrt{n}}\left[\mathbf{\widetilde V}^T\mathbcal{C}(\psi(\widetilde \alpha_k))\mathbf{\widetilde V}\right]_{kk}+o_{a.s.}(\dfrac{1}{\sqrt{n}})\right|=0, 
\end{equation*}
where 
$\left[\mathbf{\widetilde V}^T\mathbcal{C}(\psi(\widetilde \alpha_k))\mathbf{\widetilde V}\right]_{kk}$ is the $k$-th diagonal block of $\left[\mathbf{\widetilde V}^T\mathbcal{C}(\psi(\widetilde \alpha_k))\mathbf{\widetilde V}\right]$. Hence, $\chi_{k,j}$ tends to the solution of $
\left|\lambda\mathbf{I}_{n_k}-\frac1{1+cf_2(\psi(\widetilde \alpha_k))\widetilde \alpha_k}\left[\mathbf{\widetilde V}^T\mathbcal{C}(\psi(\widetilde \alpha_k))\mathbf{\widetilde V}\right]_{kk}\right|=0. 
$
As the index $j\in \mathcal{k}$ is arbitrary, it yields that almost surely, the joint distribution of $\sqrt{n}\Big\{\lambda_{n,j}-\psi(\widetilde \alpha_k), j\in \mathcal{k}\Big\}$   {converges} weakly to the joint distribution of the eigenvalues of the matrix
$	
\dfrac{1}{1+cf_2\left(\psi(\widetilde \alpha_k) \right)\widetilde \alpha_k}\left[\mathbf{\widetilde V}^T\mathbcal{C}(\psi(\widetilde \alpha_k))\mathbf{\widetilde V} \right]_{kk}.
$
Therefore, we complete the proof by returning to the original space.\hfill$\square$
\subsection{Proof of Theorem \ref{thm2}}\label{se1}
From \eqref{cal4}, we can obtain
\begin{equation}\label{cal5}
	\begin{cases}
		\mathbf{S}_{11}\boldsymbol{\beta}_{1k}+\mathbf{S}_{12}\boldsymbol{\beta}_{2k}=\lambda_{n,k}\boldsymbol{\beta}_{1k}, \\
		\mathbf{S}_{21}\boldsymbol{\beta}_{1k}+\mathbf{S}_{22}\boldsymbol{\beta}_{2k}=\lambda_{n,k}\boldsymbol{\beta}_{2k}. 
	\end{cases} 
\end{equation}
Solving the second equality of the above equation, it yields 
\begin{equation}\label{cal6}
	\boldsymbol{\beta}_{2k}=(\lambda_{n,k}\mathbf{I}_S-\mathbf{S}_{22})^{-1}\mathbf{S}_{21}\boldsymbol{\beta}_{1k}.
\end{equation}
Then, substituting  \eqref{cal6} into \eqref{cal5}, one finds 
\begin{equation*}
	[\mathbf{S}_{11}+\mathbf{S}_{12}(\lambda_{n,k}\mathbf{I}_S-\mathbf{S}_{22})^{-1}\mathbf{S}_{21}]\boldsymbol{\beta}_{1k}=\lambda_{n,k}\boldsymbol{\beta}_{1k}. 
\end{equation*}
Recalling the definition of ${\bf A}_n(\lambda)$ (see \eqref{cal7}), we get $\mathbf{A}_n(\lambda_{n,k})\boldsymbol{\beta}_{1k}=\lambda_{n,k}\boldsymbol{\beta}_{1k}$. Substituting equation \eqref{cal6} into equation \eqref{cal8} gives 
\begin{equation}\label{cal16}
	\boldsymbol{a}_k^T(\mathbf{I}_N+\mathbf{Q}_k)\boldsymbol{a}_k=\|\boldsymbol{\beta}_{1k}\|^{-2}, 
\end{equation}
where $\mathbf{Q}_k=\mathbf{S}_{12}(\lambda_{n,k}\mathbf{I}_S-\mathbf{S}_{22})^{-2}\mathbf{S}_{21}$ and $\boldsymbol{a}_k=\boldsymbol{\beta}_{1k}/\|\boldsymbol{\beta}_{1k}\|$. It is apparent from \eqref{cal10}
\begin{align}\label{cal21}
	\langle \boldsymbol{\beta}_{k},\mathbf{u}_k\rangle^2=\langle \boldsymbol{\beta}_{1k},\mathbf{\widetilde v}_k\rangle^2=\|\boldsymbol{\beta}_{1k}\|^2\cdot\langle \boldsymbol{a}_k,\mathbf{\widetilde v}_k\rangle^2.
\end{align}Subsequently, it will be shown that
\begin{align}\label{cal17}
	\boldsymbol{a}_k\xrightarrow{\text{a.s.}}\bv_k,\quad{\rm and}\quad\mathbf{Q}_k\xrightarrow{\text{a.s.}}c{f_2}(\psi(\widetilde \alpha_k))\widetilde{\boldsymbol{\Sigma}}_N.
\end{align} 

To begin with, we shall prove $\boldsymbol{a}_k\xrightarrow{\text{a.s.}}\bv_k.$ Define $\psi_n(\widetilde \alpha)$ by replacing $H$, $c$ in $\psi(\widetilde \alpha)$ with $H_n$, $c_n$ (see \eqref{cal9}). It is  known 
\begin{align*}
	\mathbf{A}_n(\lambda_{n,k})\boldsymbol{a}_k=\lambda_{n,k}\boldsymbol{a}_k,\quad \dfrac{\psi_{n}(\widetilde \alpha_k)}{\widetilde \alpha_k}\widetilde{\boldsymbol{\Sigma}}_N\bv_k=\psi_{n}(\widetilde \alpha_k)\bv_k.
\end{align*}
Using Lemma \ref{Lemma5:}, we can obtain
\begin{equation}\label{cal18}
	\boldsymbol{a}_k-\bv_k=-\mathscr{D}_{kn}\mathbf{T}_k\bv_k+\mathbf{r}_k,\ \|\mathbf{r}_k\|=O(\|\mathbf{T}_k\|^2), 
\end{equation} 
where
\begin{align*}
	&  {\mathbf{T}_k=\mathbf{A}_n(\lambda_{n,k})-\dfrac{\psi_{n}(\widetilde \alpha_k)}{\widetilde \alpha_k}\widetilde{\boldsymbol{\Sigma}}_N\in\mathbb{R}^{N\times N}, } \\
	&\mathscr{D}_{kn}=  {\sum\limits_{s\neq k}^{}\dfrac{1}{\frac{\psi_{n}(\widetilde \alpha_k)}{\widetilde \alpha_k}\widetilde \alpha_{s}-\frac{\psi_{n}(\widetilde \alpha_k)}{\widetilde \alpha_k}\widetilde \alpha_{k}}\bv_s\bv_s^T}=\dfrac{\widetilde \alpha_k}{\psi_{n}(\widetilde \alpha_k)}\sum\limits_{s\neq k}^{}\dfrac{1}{\widetilde \alpha_{s}-\widetilde \alpha_{k}}\bv_s\bv_s^T . 
\end{align*} 	
Hence, it suffices to show $\mathbf{T}_k\xrightarrow{\text{a.s.}}\mathbf{0}.$ Decompose $\mathbf{T}_k$ into the following two parts  \begin{equation*}
	\mathbf{T}_k=\mathbf{T}_{k1}+\mathbf{T}_{k2}=\bigg[\mathbf{A}_n(\psi_{n}(\widetilde \alpha_k))-\dfrac{\psi_{n}(\widetilde \alpha_k)}{\widetilde \alpha_k}\widetilde{\boldsymbol{\Sigma}}_N\bigg]+\bigg[\mathbf{A}_n(\lambda_{n,k})-\mathbf{A}_n(\psi_{n}(\widetilde \alpha_k))\bigg].
\end{equation*} 
Then, we shall prove  $\mathbf{T}_{k1}\xrightarrow{\text{a.s.}}\mathbf{0}$ and $\mathbf{T}_{k2}\xrightarrow{\text{a.s.}}\mathbf{0}$, respectively. 

To prove $\mathbf{T}_{k1}\xrightarrow{\text{a.s.}}\mathbf{0}$, we will show the entries of $\mathbf{T}_{k1}$ tend to zero almost surely. Let $\by_{Nj}$ be the transpose of $j$-th row of $\bY_N$, then the $(i,j)$-th entry of $\mathbf{T}_{k1}$ is
$$\frac1n\mathbf{y}_{Ni}^T(\mathbf{I}_{  {n}}+\mathbf{B}_n(\psi_{n}(\widetilde \alpha_k)))\mathbf{y}_{Nj}-\dfrac{\psi_{n}(\widetilde \alpha_k)}{\widetilde \alpha_k}\left(\widetilde{\boldsymbol{\Sigma}}_N\right)_{ij}.$$
Using Lemma \ref{Lemma 3:}, it yields 
\begin{equation*}
	\mathbb{E}\left|\dfrac{1}{n}\mathbf{y}_{Ni}^T(\mathbf{I}_{  {n}}+\mathbf{B}_n(\psi_{n}(\widetilde \alpha_k)))\mathbf{y}_{Nj}-\dfrac{1}{n}\left(\widetilde{\boldsymbol{\Sigma}}_N\right)_{ij} {\rm{tr}}(\mathbf{I}_{  {n}}+\mathbf{B}_n(\psi_{n}(\widetilde \alpha_k)))\right|^t \leq 
	C_tn^{-\frac{t}{2}},
\end{equation*}
which is summable for $t>2$. 
Therefore, it follows from the Borel-Cantelli lemma \begin{equation}\label{cal11}
	\dfrac{1}{n}\mathbf{y}_{Ni}^T(\mathbf{I}_{  {n}}+\mathbf{B}_n(\psi_{n}(\widetilde \alpha_k)))\mathbf{y}_{Nj}-\dfrac{1}{n}\left(\widetilde{\boldsymbol{\Sigma}}_N\right)_{ij} {\rm{tr}}(\mathbf{I}_{  {n}}+\mathbf{B}_n(\psi_{n}(\widetilde \alpha_k)))\xrightarrow{\text{a.s.}}0.
\end{equation}
Furthermore, one finds
\begin{align*}
	\dfrac{1}{n} {\rm{tr}}(\mathbf{I}_{  {n}}+\mathbf{B}_n(\psi_{n}(\widetilde \alpha_k)))=&  {\frac{1}{n}{\rm tr}\left[\psi_{n}(\widetilde \alpha_k)\left(\psi_{n}(\widetilde \alpha_k)\bI_n-\frac{1}{n}\bY_S^T\bY_S\right)^{-1}\right]}\\
	&\xrightarrow{  {\rm a.s.}}\int{\frac{\psi(\widetilde \alpha_k)}{\psi(\widetilde \alpha_k)-x}}\underline{F}_{c,H}(dx)=-\psi(\widetilde \alpha_k)\underline{m}(\psi(\widetilde \alpha_k)), 
\end{align*}
where the last equality is from the fact   {(which follows from the Silverstein equation)} that $\underline{m}(\psi(\widetilde \alpha_k))=-\dfrac{1}{\widetilde \alpha_k}$. Hence, we get
\begin{align}\label{cal12}
	\dfrac{1}{n} {\rm{tr}}(\mathbf{I}_{  {n}}+\mathbf{B}_n(\psi_{n}(\widetilde \alpha_k)))-\dfrac{\psi_{n}(\widetilde \alpha_k)}{\widetilde \alpha_k}\xrightarrow{\text{a.s.}}{0}.
\end{align}
Combining with \eqref{cal11} and \eqref{cal12}, we deduce $$\frac1n\mathbf{y}_{Ni}^T(\mathbf{I}_{  {n}}+\mathbf{B}_n(\psi_{n}(\widetilde \alpha_k)))\mathbf{y}_{Nj}-\dfrac{\psi_{n}(\widetilde \alpha_k)}{\widetilde \alpha_k}\left(\widetilde{\boldsymbol{\Sigma}}_N\right)_{ij}\xrightarrow{\text{a.s.}}0.$$

For $\mathbf{T}_{k2}$, we give some notations firstly. Let
\begin{equation*}
	\mathbf{G}_n(\lambda)=\left(\lambda \mathbf{I}_n-\underline{\bS}_{22}\right)^{-1} ,\qquad \mathbf{A}_n(\lambda)=\frac{\lambda}{n}\mathbf{Y}_N \mathbf{G}_n(\lambda)\mathbf{Y}_N^T ,
\end{equation*}
where $\mathbf{I}_n+\mathbf{B}_n(\lambda)=\mathbf{I}_n+\underline{\bS}_{22}\mathbf{G}_n(\lambda)=\lambda \mathbf{G}_n(\lambda)$. 
By the formula $$\mathbf{A}^{-1}-\mathbf{B}^{-1}=\mathbf{A}^{-1}(\mathbf{B}-\mathbf{A})\mathbf{B}^{-1},$$ we get
\begin{equation}\label{cal13}
	\mathbf{G}_n(\lambda)-\mathbf{G}_n(\psi)=-(\lambda-\psi)\mathbf{G}_n(\lambda)\mathbf{G}_n(\psi), 
\end{equation}
\begin{equation*}
	\lambda \mathbf{G}_n(\lambda)-\psi \mathbf{G}_n(\psi)=-(\lambda-\psi)\underline{\bS}_{22}\mathbf{G}_n(\lambda)\mathbf{G}_n(\psi). 
\end{equation*}
Hence, 
\begin{align}\label{cal15}
	\mathbf{T}_{k2}&= -(\lambda_{n,k}-\psi_{n}(\widetilde \alpha_k))\dfrac{1}{n}\mathbf{Y}_N\underline{\bS}_{22}\mathbf{G}_n(\lambda_{n,k})\mathbf{G}_n(\psi_{n}(\widetilde \alpha_k))\mathbf{Y}_N^T \notag\\
	&=-(\lambda_{n,k}-\psi_{n}(\widetilde \alpha_k))\dfrac{1}{n}\mathbf{Y}_N\underline{\bS}_{22}\mathbf{G}_n^2(\psi_{n}(\widetilde \alpha_k))\mathbf{Y}_N^T\\
	&\quad+(\lambda_{n,k}-\psi_{n}(\widetilde \alpha_k))^2\dfrac{1}{n}\mathbf{Y}_N\underline{\bS}_{22}\mathbf{G}_n(\lambda_{n,k})\mathbf{G}_n^2(\psi_{n}(\widetilde \alpha_k))\mathbf{Y}_N^T.\notag
\end{align}
By   {Proposition \ref{pro1}, } we have for any $\epsilon>0$, 
\begin{align}\label{cal14}
	\lambda_{n,k}>\psi(\widetilde \alpha_k)-\epsilon,\quad a.s..
\end{align}It follows that from \eqref{cal13} and \eqref{cal14}
$\mathbf{G}_n(\lambda_{n,k}) < \mathbf{G}_n(\psi(\widetilde \alpha_k)-\epsilon),$ which implies  
\begin{equation*}
	\dfrac{1}{n}\mathbf{Y}_N\mathbf{
		J}_{nst}(\lambda_{n,k},\psi_{n}(\widetilde \alpha_k))\mathbf{Y}_N^T < \dfrac{1}{n}\mathbf{Y}_N\mathbf{J}_{nst}(\psi(\widetilde \alpha_k)-\epsilon , \psi_{n}(\widetilde \alpha_k))\mathbf{Y}_N^T, 
\end{equation*}
where 
$
\mathbf{J}_{nst}(\lambda , \psi)=\underline{\bS}_{22}\mathbf{G}_n^s(\lambda)\mathbf{G}_n^t(\psi). 
$
Similar to \eqref{cal11}, we have
\begin{align*}
	\dfrac{1}{n}\mathbf{Y}_N\mathbf{J}_{nst}(\psi(\widetilde \alpha_k)-\epsilon , \psi_{n}(\widetilde \alpha_k))\mathbf{Y}_N^T-\dfrac{1}{n}\widetilde{\boldsymbol{\Sigma}}_N {\rm{tr}}\bigg(\mathbf{J}_{nst}(\psi(\widetilde \alpha_k)-\epsilon , \psi_{n}(\widetilde \alpha_k))\bigg)\xrightarrow{\text{a.s.}}0.
\end{align*}
Considering the spectral norm of $\mathbf{J}_{nst}(\psi(\widetilde \alpha_k)-\epsilon , \psi_{n}(\widetilde \alpha_k))$ is almost surely bounded, it yields
\begin{equation*}
	\dfrac{1}{n}\mathbf{Y}_N\mathbf{J}_{nst}(\lambda_{n,k},\psi_{n}(\widetilde \alpha_k))\mathbf{Y}_N^T=c\int\frac{x}{(\psi_{n}(\widetilde \alpha_k)-x)^{s+t}}F_{c,H}(dx)\widetilde{\boldsymbol{\Sigma}}_N+o_{a.s.}(1). 
\end{equation*}
By Theorem \ref{thm1}, \eqref{cal15}, and the above equality, we conclude
\begin{equation}\label{cal20}
	\mathbf{T}_{k2}=-(\lambda_{n,k}-\psi_{n}(\widetilde \alpha_k))\cdot \left[cf_2(\psi_{n}(\widetilde \alpha_k))\widetilde{\boldsymbol{\Sigma}}_N+o_{a.s.}(1)\right]\xrightarrow{\text{a.s.}}\mathbf{0}. 
\end{equation}

Now, we are in position to find the almost surely limit of $\mathbf{Q}_k$. Let $\mathbf{L}_S(\lambda)=(\lambda \mathbf{I}_S-\mathbf{S}_{22})^{-1}$, then we have the following decomposition
\begin{align*}\mathbf{Q}_k=\mathbf{S}_{12}\mathbf{L}_S^2(\lambda_{n,k})\mathbf{S}_{21}&=\mathbf{S}_{12}\mathbf{L}_S^2(\psi(\widetilde \alpha_k))\mathbf{S}_{21}+\mathbf{S}_{12}[\mathbf{L}_S^2(\lambda_{n,k})-\mathbf{L}_S^2(\psi(\widetilde \alpha_k))]\mathbf{S}_{21}\\
	&=\dfrac{1}{n}\mathbf{Y}_N\widetilde{\mathbf{L}}_{n1}\mathbf{Y}_N^T-(\lambda_{n,k}-\psi(\widetilde \alpha_k))\dfrac{1}{n}\mathbf{Y}_N\widetilde{\mathbf{L}}_{n2}\mathbf{Y}_N^T
	=\mathbf{Q}_{k1}+\mathbf{Q}_{k2}, 
\end{align*}
where 
\begin{equation*}
	\widetilde{\mathbf{L}}_{n1}=\dfrac{1}{n}\mathbf{Y}_S^T\mathbf{L}_S^2(\psi(\widetilde \alpha_k))\mathbf{Y}_S , \quad \widetilde{\mathbf{L}}_{n2}=\dfrac{1}{n}\mathbf{Y}_S^T\left[(\lambda_{n,k}+\psi(\widetilde \alpha_k))\mathbf{I}_S-2\mathbf{S}_{22}\right]\mathbf{L}_S^2(\lambda_{n,k})\mathbf{L}_S^2(\psi(\widetilde \alpha_k))\mathbf{Y}_S
	.
\end{equation*}
It can be verified $\widetilde{\mathbf{L}}_{n1}$ and $\widetilde{\mathbf{L}}_{n2}$ are bounded almost surely,   {since $\lambda_{n,k}\xrightarrow{\text{a.s.}}\psi(\widetilde{\alpha}_k)\notin\Gamma_{F_{c,H}}$, and almost surely, no eigenvalue of $\bS_{22}$ lies outside $\Gamma_{F_{c,H}}$}.  
Using the same procedure   {to that leading to} \eqref{cal11},   {we have}
$
\dfrac{1}{n}\mathbf{Y}_N\widetilde{\mathbf{L}}_{n1}\mathbf{Y}_N^T-\dfrac{1}{n}\widetilde{\boldsymbol{\Sigma}}_N {\rm{tr}}\bigg(\widetilde{\mathbf{L}}_{n1}\bigg)\xrightarrow{\text{a.s.}}0
$
and
$
\dfrac{1}{n}\mathbf{Y}_N\widetilde{\mathbf{L}}_{n2}\mathbf{Y}_N^T-\dfrac{1}{n}\widetilde{\boldsymbol{\Sigma}}_N {\rm{tr}}\bigg(\widetilde{\mathbf{L}}_{n2}\bigg)\xrightarrow{\text{a.s.}}0.
$
Applying   {Proposition} \ref{pro1}, we find $\mathbf{Q}_{k2}  {=-(\lambda_{n,k}-\psi(\widetilde \alpha_k))\dfrac{1}{n}\widetilde{\boldsymbol{\Sigma}}_N {\rm{tr}}\bigg(\widetilde{\mathbf{L}}_{n2}\bigg)+o_{a.s.}(1)}\xrightarrow{\text{a.s.}}0$. Moreover, one gets
$
\mathbf{Q}_{k1}=\dfrac{{\rm{tr}}(\widetilde{\mathbf{L}}_{n1})}{n}\widetilde{\boldsymbol{\Sigma}}_N+o_{a.s.}(1)\xrightarrow{a.s.}c{f_2}(\psi(\widetilde \alpha_k))\widetilde{\boldsymbol{\Sigma}}_N . 
$
Hence, we conclude that
$
\mathbf{Q}_{k}\xrightarrow{a.s.}c{f_2}(\psi(\widetilde \alpha_k))\widetilde{\boldsymbol{\Sigma}}_N . 
$ Furthermore, combining with \eqref{cal16}-\eqref{cal17}, we can obtain $$\langle\boldsymbol{\beta}_{k},{\bf u}_k\rangle^2=\frac{\langle \boldsymbol{a}_k,\mathbf{\widetilde v}_k\rangle^2}{1+\boldsymbol{a}_k^T\mathbf{Q}_k\boldsymbol{a}_k}
\xrightarrow{\text{a.s.}}\dfrac{1}{1+c{f_2}(\psi(\widetilde \alpha_k))\widetilde \alpha_{k}}. $$

  {Specifically, when $\psi'(\widetilde{\alpha}_k)\leq0$, we aim to prove $\langle \boldsymbol{\beta}_{k},\mathbf{u}_k\rangle^2=\|\boldsymbol{\beta}_{1k}\|^2\cdot\langle \boldsymbol{a}_k,\mathbf{\widetilde v}_k\rangle^2 \xrightarrow{\text{a.s.}}
	0$, and by \eqref{cal16}, it suffices to prove $\boldsymbol{a}_k^T\mathbf{Q}_k\boldsymbol{a}_k\xrightarrow{\text{a.s.}}\infty$. We obtain this result by proving $\lambda_{\rm min}(\mathbf{Q}_k)\xrightarrow{\text{a.s.}}\infty$. We give the proof when  {  $\widetilde{\alpha}_k > \lambda_{\rm max}(\widetilde{\boldsymbol{\Sigma}}_S)$} in what follows. The case {  $\widetilde{\alpha}_k <\lambda_{\rm min}(\widetilde{\boldsymbol{\Sigma}}_S)$} can be proved similarly.  Denote $\mathbf{Q}_{k\epsilon}(\lambda_{n,k})=\mathbf{S}_{12}\left[(\lambda_{n,k}\mathbf{I}_S-\mathbf{S}_{22})^{2}+\epsilon^2\bI_S\right]^{-1}\mathbf{S}_{21}$, for $\epsilon>0$.  Combining Proposition \ref{pro1}, we perform the following decomposition \begin{align*}
		\liminf\lambda_{\rm min}(\mathbf{Q}_k)\geq\liminf\lambda_{\rm min}(\mathbf{Q}_{k\epsilon}(\lambda_{n,k}))=\liminf\lambda_{\rm min}(\mathbf{Q}_{k\epsilon}(\sup \Gamma_{F_{c,H}})+\Delta_{k\epsilon}), 
	\end{align*} where $\Delta_{k\epsilon}\triangleq\mathbf{Q}_{k\epsilon}(\lambda_{n,k})-\mathbf{Q}_{k\epsilon}(\sup \Gamma_{F_{c,H}})$. Using the same procedure to that leading to \eqref{cal11}, we have \begin{align*}
		\mathbf{Q}_{k\epsilon}(\sup \Gamma_{F_{c,H}})\xrightarrow{\text{a.s.}}c\displaystyle \int \dfrac{x}{(\sup \Gamma_{F_{c,H}} - x)^2+\epsilon^2} \, F_{c,H}(dx)\widetilde{\boldsymbol{\Sigma}}_N\triangleq c{f_2}(\sup \Gamma_{F_{c,H}},\epsilon)\widetilde{\boldsymbol{\Sigma}}_N. 
	\end{align*}
	Furthermore, we derive $\Delta_{k\epsilon}\xrightarrow{\text{a.s.}}0$ in the same manner as the proof of $\mathbf{Q}_{k2}\xrightarrow{\text{a.s.}}0$. Since $\lambda_{\rm min}(\cdot)$ is a continuous function on $N\times N$ matrices, we obtain that \begin{align*}
		\liminf\lambda_{\rm min}(\mathbf{Q}_k)\geq c{f_2}(\sup \Gamma_{F_{c,H}},\epsilon)\liminf\lambda_{\rm min}(\widetilde{\boldsymbol{\Sigma}}_N), 
	\end{align*}
	and since ${f_2}(\sup \Gamma_{F_{c,H}},\epsilon)\geq{f_2}(\sup \Gamma_{F_{c,H}}+\epsilon)$ and ${f_2}(\sup \Gamma_{F_{c,H}}+\epsilon)\rightarrow\infty$ as $\epsilon\rightarrow0$, we get $\lambda_{\rm min}(\mathbf{Q}_k)\xrightarrow{\text{a.s.}}\infty$. }\hfill$\square$
\subsection{Proof of Theorem \ref{thm3}}
By \eqref{cal18}, it follows that
\begin{align}\label{cal22}
	\sqrt{n}(\boldsymbol{a}_k-\bv_k)&=-\sqrt{n}\mathscr{D}_{kn}\mathbf{T}_k\bv_k+\sqrt{n}\boldsymbol{r}_k, \ \|\mathbf{r}_k\|=O(\|\mathbf{T}_k\|^2),   
\end{align}
  {where 
	\begin{equation*}
		\mathbf{T}_k=\mathbf{T}_{k1}+\mathbf{T}_{k2}=\bigg[\mathbf{A}_n(\psi_{n}(\widetilde \alpha_k))-\dfrac{\psi_{n}(\widetilde \alpha_k)}{\widetilde \alpha_k}\widetilde{\boldsymbol{\Sigma}}_N\bigg]+\bigg[\mathbf{A}_n(\lambda_{n,k})-\mathbf{A}_n(\psi_{n}(\widetilde \alpha_k))\bigg].
\end{equation*}}
For $\mathbf{T}_{k1}$, we have
  {\begin{align*}
		\sqrt{n}\mathbf{T}_{k1}&=\sqrt{n}\dfrac{\psi_{n}(\widetilde \alpha_{k})}{n}\left[\mathbf{Y}_N\mathbf{G}_n(\psi_{n}(\widetilde \alpha_k))\mathbf{Y}_N^T-{\rm{tr}}\bigg(\mathbf{G}_n(\psi_{n}(\widetilde \alpha_k))\bigg)\widetilde{\boldsymbol{\Sigma}}_N\right]\\&+\sqrt{n}\dfrac{\psi_{n}(\widetilde \alpha_k)}{n}\left[{\rm{tr}}\bigg(\mathbf{G}_n(\psi_{n}(\widetilde \alpha_k))\bigg)+n\cdot \underline m_n(\psi_{n}(\widetilde \alpha_k))\right]\widetilde{\boldsymbol{\Sigma}}_N\\
		&=\mathbf{C}_n(\psi_{n}(\widetilde \alpha_k))+o_p(1)
		, \end{align*}}
where the last equality   {follows from Theorem 2.3 in} \cite{li2024spectral}. 
It follows that from \eqref{cal20} and Theorem \ref{thm2}
  {\begin{equation*}
		\sqrt{n}\mathbf{T}_{k2}
		=-\sqrt{n}(\lambda_{n,k}-\psi_{n}(\widetilde \alpha_k))\cdot cf_2(\psi_{n}(\widetilde \alpha_k))\widetilde{\boldsymbol{\Sigma}}_N+o_p(1). 
\end{equation*}} 
Hence, 	we get   {\begin{align*}
		\sqrt{n}\mathbf{T}_{k}
		&=\mathbf{C}_n(\psi_{n}(\widetilde \alpha_k))-\sqrt{n}(\lambda_{n,k}-\psi_{n}(\widetilde \alpha_k))\cdot cf_2(\psi_{n}(\widetilde \alpha_k))\widetilde{\boldsymbol{\Sigma}}_N+o_p(1)
		.\end{align*}}
  {Since
	\begin{align*}
		\mathscr{D}_{kn} \widetilde{\boldsymbol{\Sigma}}_N \bv_k&=\dfrac{\widetilde \alpha_k}{\psi_{n}(\widetilde \alpha_k)}\sum\limits_{s\neq k}^{}\dfrac{1}{\widetilde \alpha_{s}-\widetilde \alpha_{k}}\bv_s\bv_s^T\widetilde{\boldsymbol{\Sigma}}_N \bv_k=\dfrac{\widetilde \alpha^2_k}{\psi_{n}(\widetilde \alpha_k)}\sum\limits_{s\neq k}^{}\dfrac{1}{\widetilde \alpha_{s}-\widetilde \alpha_{k}}\bv_s\bv_s^T\bv_k,
	\end{align*} 
	and given that when \(s \neq k\), \(\bv_s\) and \(\bv_k\) are orthogonal, we have \(\mathscr{D}_{kn} \widetilde{\boldsymbol{\Sigma}}_N \bv_k = 0\).}
  {We then deduce from \eqref{cal22} and the above equality that}
\begin{align*}
	\sqrt{n}(\boldsymbol{a}_k-\bv_k)=-\mathscr{D}_{kn}\mathbf{C}_n(\psi_{n}(\widetilde \alpha_k))\bv_{k}+o_p(1).
\end{align*}
It is known that
\begin{equation*}
	\mathbf{C}_n(\psi_{n}(\widetilde \alpha_k))\xrightarrow{\text{d}}\mathbcal{C}(\psi(\widetilde \alpha_k)),\qquad \mathscr{D}_{kn}\rightarrow\mathscr{D}_{k}.
\end{equation*}
Therefore, by the Slutsky's theorem, we get
\begin{equation*}
	\sqrt{n}(\boldsymbol{a}_k-\bv_k) \xrightarrow{\text{d}}N(\boldsymbol{0},\boldsymbol{\Gamma}_k), 
\end{equation*}
where
\begin{equation*}
	\boldsymbol{\Gamma}_k=\mathscr{D}_{k}\text{Cov}(\mathbcal{C}(\psi(\widetilde \alpha_k))\bv_k)\mathscr{D}_{k}. 
\end{equation*}
  {It can be checked that}
\begin{align*}
	\boldsymbol{\Gamma}_k=\dfrac{\widetilde \alpha_k^2}{\psi^2(\widetilde \alpha_k)}\sum\limits_{f_1\neq k}\sum\limits_{f_2\neq k}^{}\dfrac{1}{(\widetilde \alpha_{f_1}-\widetilde \alpha_{k})(\widetilde \alpha_{f_2}-\widetilde \alpha_{k})}\bv_{f_1}\bv_{f_1}^T \text{Cov}(\mathbcal{C}(\psi(\widetilde \alpha_k))\bv_k)\bv_{f_2}\bv_{f_2}^T. 
\end{align*} 

By Theorem \ref{thm2}, we see
\begin{align*}
	&\bv_{f_1}^T \text{Cov}(\mathbcal{C}(\psi(\widetilde \alpha_k))\bv_k)\bv_{f_2}\\
	=&\tau_0\left(\psi(\widetilde \alpha_k) \right)\sum_{j,\ell,s,t}\widetilde{v}_{f_1j}\widetilde{v}_{k\ell}\widetilde{v}_{ks}\widetilde{v}_{f_2t}\mathbb{E}\bigg(y_{j1}y_{s1}y_{l1}y_{t1}\bigg)\\
	&-\tau_0\left(\psi(\widetilde \alpha_k) \right)\bv_{f_1}^T\widetilde{\boldsymbol{\Sigma}}_N\bv_{k}\bv_{k}^T\widetilde{\boldsymbol{\Sigma}}_N\bv_{f_2} \\
	&+\Big(\nu_0\left(\psi(\widetilde \alpha_k) \right)-\tau_0\left(\psi(\widetilde \alpha_k) \right)\Big)\bigg[\bv_{f_1}^T\widetilde{\boldsymbol{\Sigma}}_N\bv_{f_2}\bv_{k}^T\widetilde{\boldsymbol{\Sigma}}_N\bv_{k}+\bv_{f_1}^T\widetilde{\boldsymbol{\Sigma}}_N\bv_{k}\bv_{k}^T\widetilde{\boldsymbol{\Sigma}}_N\bv_{f_2}\bigg] \\
	=&\tau_0\left(\psi(\widetilde \alpha_k) \right)\bv_{f_1}^T\mathscr{C}^k\bv_{f_2}+\widetilde \alpha_k\widetilde \alpha_{f_2}\Big(\nu_0\left(\psi(\widetilde \alpha_k) \right)-\tau_0\left(\psi(\widetilde \alpha_k) \right)\Big)\bv_{f_1}^T\bv_{f_2}.
\end{align*}
Hence, it implies
\begin{align*}
	\boldsymbol{\Gamma}_k
	=&\tau_0\left(\psi(\widetilde \alpha_k) \right)\mathscr{D}_{k}\mathscr{C}^k\mathscr{D}_{k}+\dfrac{\Big(\nu_0\left(\psi(\widetilde \alpha_k) \right)-\tau_0\left(\psi(\widetilde \alpha_k) \right)\Big)}{\psi^2(\widetilde \alpha_k)}\sum\limits_{f\neq k}\dfrac{\widetilde \alpha_k^3\widetilde \alpha_{f}}{(\widetilde \alpha_{f}-\widetilde \alpha_{k})^2}\bv_{f}\bv_{f}^T. 
\end{align*}
Using the Delta method, we get
\begin{equation*}
	\sqrt{n}({\bV^T}\boldsymbol{a}_k-\be_k) \xrightarrow{\text{d}}N(\boldsymbol{0},\bV^T\boldsymbol{\Gamma}_k\bV), 
\end{equation*}
where
\begin{align*}
	\bV^T\boldsymbol{\Gamma}_k\bV
	=&\tau_0\left(\psi(\widetilde \alpha_k) \right)\bV^T\mathscr{D}_{k}\mathscr{C}^k\mathscr{D}_{k}\bV+\dfrac{\Big(\nu_0\left(\psi(\widetilde \alpha_k) \right)-\tau_0\left(\psi(\widetilde \alpha_k) \right)\Big)}{\psi^2(\widetilde \alpha_k)}\sum\limits_{f\neq k}\dfrac{\widetilde \alpha_k^3\widetilde \alpha_{f}}{(\widetilde \alpha_{f}-\widetilde \alpha_{k})^2}\be_{f}\be_{f}^T. 
	\tag*{$\square$}
\end{align*}

\section{Some lemmas}\label{sec5}
This section contains several lemmas used in the proof. The first lemma is crucial for determining the limiting distribution of spiked eigenvalues and eigenvectors.
\begin{lemma}\label{Lemma 1: }
	Let $\mathbf{M} = (m_{ij})$ be an $n \times n$ non-random symmetric matrix. Additionally, consider a $2g \times n$ data matrix \(\left(\Xi^T, \mathbcal{H}^T\right)^T\) sampled from a $2g$-dimensional incomplete population, where
	\begin{equation*}
		\Xi^T =
		\left(d_{ij}x_{ij}\right)_{n\times g}
		= \left(\boldsymbol{\xi}_1 , \dots ,\boldsymbol{\xi}_g\right), \ 
		\mathbcal{H}^T =
		\left(q_{ij}y_{ij}\right)_{n\times g}
		= \left(\boldsymbol{\eta}_1 ,\dots ,\boldsymbol{\eta}_g\right), 
	\end{equation*}
	and $d_{ij}\sim{B}(1,\theta^{\xi}_i)$, $(i=1, \ldots, g, j=1, \ldots, n)$, $q_{ij}\sim {B}(1, \theta^{\eta}_{i})$, $(i=1, \ldots, g, j=1, \ldots, n)$. Let $\rho_i=\mathbb{E}(d_{i1}q_{i1}x_{i1}y_{i1})$. Assume that the $2g$-dimensional complete population has an expectation of $\mathbf{0}$ and a bounded fourth moment, where $g$ is finite, and that the following limits exist:
	\begin{align*}
		&\tau=\lim_{n \to \infty}\dfrac{1}{n}\sum_{i=1}^{n}m_{ii}^2, \ \
		\nu=\lim_{n \to \infty}\dfrac{1}{n}{\rm tr}(\mathbf{M}^2)
		.\end{align*}
	Let the $g$-dimensional real random vector $\mathbf{z}^T=\left(z_1,\dots,z_g \right)$
	where $z_l=\dfrac{1}{\sqrt{n}}(\boldsymbol{\xi}_l^T \mathbf{M} \boldsymbol{\eta}_l- \rho_l{\rm{tr}}\mathbf{M})$, $(l=1, \ldots, g )$. Then, we have 
	\[\mathbf{z}\xrightarrow{\text{d}}\mathbcal{g}\sim N_g(\mathbf{0},\mathbf{G}), \]  where  $\mathbf{G}=\mathbf{G}_1+\mathbf{G}_2+\mathbf{G}_3$, 
	\begin{align*}
		\left(\mathbf{G}_1\right)_{ij}
		&=\tau\{\mathbb{E}\left(d_{i1}q_{i1}d_{j1}q_{j1}x_{i1}y_{i1}x_{j1}y_{j1} \right)-\mathbb{E}\left(d_{i1}x_{i1}q_{i1}y_{i1}\right)\mathbb{E}\left(d_{j1}x_{j1}q_{j1}y_{j1} \right)\}, \\
		\left(\mathbf{G}_2\right)_{ij}
		&=(\nu-\tau)\{\mathbb{E}\left(d_{i1}x_{i1}q_{j1}y_{j1}\right)\mathbb{E}\left(d_{j1}x_{j1}q_{i1}y_{i1}\right) \}, \\
		\left(\mathbf{G}_3\right)_{ij}
		&=(\nu-\tau)\{\mathbb{E}\left(d_{i1}x_{i1}d_{j1}x_{j1}\right)\mathbb{E}\left(q_{i1}y_{j1}q_{j1}y_{j1}\right) \}.
	\end{align*}
\end{lemma}
\begin{Proof}
	  {The detailed proof of this lemma is provided in the supplementary material \cite{cheng2026supplement}.}
\end{Proof}

\begin{lemma}\label{Lemma 2:}
Under the conditions of Theorem \ref{thm1},	we have for $  {k}\ge1$
\begin{align*}
	\dfrac{1}{n}{\rm tr}\Big[\mathbf{B}_n^{  {k}}(\lambda)\Big] \xrightarrow{\text{p}} c \displaystyle \int \dfrac{x^{  {k}}}{(\lambda - x)^{  {k}}} \, F_{c,H}(dx) , \quad
	\dfrac{1}{n}\sum_{i=1}^{n}(\mathbf{B}_n(\lambda))_{ii}^2 \xrightarrow{\text{p}} \left(\int \dfrac{  {x}}{\lambda - x} \, \underline{F}_{c,H}(dx) \right)^2 
	, 
\end{align*}
where
\begin{align*}
	  {\mathbf{B}_n(\lambda)=\dfrac{1}{n}\mathbf{Y}_S^T(\lambda \mathbf{I}_S-\frac{1}{n}\mathbf{Y}_S\mathbf{Y}_S^T)^{-1}\mathbf{Y}_S.}
\end{align*}
\end{lemma}
\begin{Proof}
  {Since $N$ is a fixed number, it follows from Assumption (a), $p/n=(N+S)/n\rightarrow c$ implies that $S/n\rightarrow c$. }
Let $F_{n2}$ denote the ESD of $\bS_{22}$, and then it follows from   {Theorem 2.1 of } \cite{li2024spectral}
\begin{align*}
	\dfrac{1}{n}{\rm tr}\Big[\mathbf{B}_n^{  {k}}(\lambda)\Big]
	&=\dfrac{S}{n} \int \dfrac{x^{  {k}}}{(\lambda - x)^{  {k}}} \, F_{n2}(dx) 
	\xrightarrow{\text{p}}c  \int \dfrac{x^{  {k}}}{(\lambda - x)^{  {k}}} \, F_{c,H}(dx). 
\end{align*}
  {We recall the resolvent identity
	\[
	(\lambda\bI - \mathbf{B}\mathbf{A})^{-1} = \lambda^{-1}\left(\bI + \mathbf{B}(\lambda\bI - \mathbf{A}\mathbf{B})^{-1}\mathbf{A}\right),
	\]
	which implies
	\[
	\mathbf{A}(\lambda\bI - \mathbf{B}\mathbf{A})^{-1}\mathbf{B} = (\lambda\bI - \mathbf{A}\mathbf{B})^{-1}\mathbf{A}\mathbf{B},
	\]
	using $\mathbf{M}(\lambda\bI - \mathbf{M})^{-1} = \lambda(\lambda\bI - \mathbf{M})^{-1} - \bI$. With $\underline{\bS}_{22} = \frac{1}{n}\bY_S^T \bY_S$, we have $\mathbf{B}_n(\lambda) = (\lambda \bI_n - \underline{\bS}_{22})^{-1}\underline{\bS}_{22}$, so that $\bI_n + \mathbf{B}_n(\lambda) = \lambda(\lambda \bI_n - \underline{\bS}_{22})^{-1}$.
	We obtain}
\begin{align*}
	(\mathbf{I}_n+\mathbf{B}_n(\lambda))_{ii}&=\be_i^T(\mathbf{I}_n+\mathbf{B}_n(\lambda))\be_i=  {\be_i^T\lambda\left(\lambda\mathbf{I}_n-\underline{\bS}_{22}\right)^{-1}\underline{\bS}_{22}\be_i}
	\xrightarrow{p}  \int \dfrac{\lambda}{\lambda - x} \, \underline{F}_{c,H}(dx). 
\end{align*}
Hence, we get
\begin{equation*}
\dfrac{1}{n}\sum_{i=1}^{n}(\mathbf{B}_n(\lambda))_{ii}^2 \xrightarrow{\text{p}} \left(\int \dfrac{  {x}}{\lambda - x} \, \underline{F}_{c,H}(dx)\right)^2.\qedhere
\end{equation*}
\end{Proof}
\begin{lemma}\label{Lemma 3:}
Assume the independent zero-mean vectors
\begin{equation*}
\begin{pmatrix}
	x_i \\
	y_i
\end{pmatrix},\qquad i=1,\dots,n,
\end{equation*}
have covariance matrix
\begin{equation*}
\operatorname{Cov}\left(\begin{pmatrix}
	x_i \\
	y_i
\end{pmatrix}\right)=\begin{pmatrix}
	1 & \rho \\
	\rho & 1
\end{pmatrix},\qquad \max\left\{\mathbb{E}|x_i|^p,\mathbb{E}|y_i|^p\right\}\leq v_p,\quad {\rm for}\ p\geq2.
\end{equation*}
Let  $\mathbf{B}$ be an $n$-dimensional symmetric matrix, $\mathbf{x}=(x_1,\dots,x_n)^T$, $\mathbf{y}=(y_1,\dots,y_n)^T$, then
\begin{equation*}
	\mathbb{E}\left|\mathbf{x}^T\mathbf{B}\mathbf{y}-\rho\cdot {\rm{tr}}(\mathbf{B})\right|^p \leq C_p\left[\bigg(v_4{\rm{tr}}(\mathbf{B}^2)\bigg)^{\frac{p}{2}}+v_{2p}{\rm{tr}}(\mathbf{B}^p)\right].
\end{equation*}
\end{lemma}
\begin{Proof}
Write $\mathbf{x}^T\mathbf{B}\mathbf{y}=(1+\rho)\mathbf{z}^T\mathbf{B}\mathbf{z}-\frac{1}{2}\mathbf{x}^T\mathbf{B}\mathbf{x}-\frac{1}{2}\mathbf{y}^T\mathbf{B}\mathbf{y}$
where $\mathbf{z}=\dfrac{(\mathbf{x}+\mathbf{y})}{\sqrt{2+2\rho}}$.
By Lemma 2.2 in \cite{2004CLT} and $(a+b+c)^p\leq 3^{p-1}(a^p+b^p+c^p)$, we can get \begin{align*}
	\mathbb{E}\left|\mathbf{x}^T\mathbf{B}\mathbf{y}-\rho {\rm{tr}}\mathbf{B}\right|^p\le&   {3^{p-1}(1 + \rho^p)} \mathbb{E}\left|\mathbf{z}^T\mathbf{B}\mathbf{z}- {\rm{tr}}\mathbf{B}\right|^p+  {(3^{p-1}/2^p)}\mathbb{E}\left|\mathbf{x}^T\mathbf{B}\mathbf{x}- {\rm{tr}}\mathbf{B}\right|^p\\
	&+  {(3^{p-1}/2^p)}\mathbb{E}\left|\mathbf{y}^T\mathbf{B}\mathbf{y}- {\rm{tr}}\mathbf{B}\right|^p\\
	\le&C_p\left[(v_4{\rm{tr}}\mathbf{B}^2)^{\frac{p}{2}}+v_{2p}{\rm{tr}}\mathbf{B}^p\right], 
\end{align*}
  {where the constant $C_p$ depends on $\rho$. }
\end{Proof} 

\begin{lemma}\label{Lemma5:}(Lemma 13 in \cite{2018Notes}) Let $\mathbf{A}$ and $\mathbf{B}$ are $N$-dimensional symmetric matrices. The eigenvalues of $\mathbf{A}$ are $\lambda_1(\mathbf{A}) \geq \dots \geq \lambda_N(\mathbf{A})$ and the corresponding eigenvectors are $\mathbf{v}_1(\mathbf{A}) , \dots , \mathbf{v}_N(\mathbf{A})$. Define
\begin{equation*}
	\mathscr{D}_k(\mathbf{A})=\sum\limits_{s\neq k}^{}\dfrac{1}{\lambda_s(\mathbf{A})-\lambda_k(\mathbf{A})}\mathbf{V}_s(\mathbf{A}) , 
\end{equation*} 
where $\mathbf{ V}_s(\mathbf{A})=\mathbf{v}_s(\mathbf{A})\mathbf{v}_s(\mathbf{A})^T$. When $\mathscr{D}_k(\mathbf{A})$ satisfies
$
\mathscr{D}_k(\mathbf{A})(\mathbf{A}-\lambda_k(\mathbf{A})\mathbf{I}_N)=\mathbf{I}_N-\mathbf{V}_k(\mathbf{A}), 
$ we call $\mathscr{D}_k(\mathbf{A})$ is the resolvent of $\mathbf{A}$ "evaluated at $\lambda_k(\mathbf{A})$". Let
$
\triangle_{k}(\mathbf{A})={\rm min}\left\{\left|\lambda_i(\mathbf{A})-\lambda_k(\mathbf{A})\right|:1\leq i \neq k \leq N\right\}.
$
If $\|\mathbf{B}\|<\triangle_{k}(\mathbf{A})/3$, the $\lambda_k(\mathbf{A}+\mathbf{B})$ is also simple and the selection of symbols is to satisfy $\delta_k=\mathbf{v}_k(\mathbf{A})^T\mathbf{v}_k(\mathbf{A}+\mathbf{B})\geq0$. Then we have 
\begin{equation*}
	\mathbf{v}_k(\mathbf{A}+\mathbf{B})-\mathbf{v}_k(\mathbf{A})=-\mathscr{D}_k(\mathbf{A})\mathbf{B}\mathbf{v}_k(\mathbf{A})+\mathbf{r}_k \quad and \quad \|\mathbf{r}_k\|\leq 10\triangle_{k}^{-2}(\mathbf{A})\|\mathbf{B}\|^{2}.
\end{equation*}
\end{lemma}
\begin{acks}
All authors contributed equally and are listed in alphabetical order by surname. Huiqin Li and Yanqing Yin are the co-corresponding authors.

Huiqin Li is supported by the Open Project of the Nanjing Audit University Joint Lab for Statistics and Finance (Grant No. 2026JLSF303) and the Jiangsu Provincial Key Discipline Construction Project (Statistics). Yanqing Yin is partially supported by the National Natural Science Foundation of China (NSFC) under Grant No. 12271065. Zhixiang Zhang was partially supported by FDCT/0050/2025/ITP1, NSFC Grant No. 12401331, and the University of Macau under Grant No. MYRG-GRG2024-00260-FST-UMDF.
\end{acks}

\begin{supplement}
\stitle{Supplement to ``Limiting eigen-structure of spiked sample covariance matrices under missing observations''}
\sdescription{Simulation results for the empirical sizes and powers of the test statistic $\hat{\mathbcal{T}}_N$ under different missing probabilities, together with the detailed proof of Lemma \ref{Lemma 1: }.}
\end{supplement}

\bibliographystyle{imsart-number}
\bibliography{references}

\begin{thebibliography}{17}

\bibitem{1958An}
\begin{bbook}[author]
\bauthor{\bsnm{Anderson},~\bfnm{T.~W.}\binits{T.~W.}}
(\byear{1958}).
\btitle{An Introduction to Multivariate Statistical Analysis}.
\bpublisher{Wiley}, \baddress{New York}.
\end{bbook}
\endbibitem

\bibitem{2004CLT}
\begin{barticle}[author]
\bauthor{\bsnm{Bai},~\bfnm{Z.~D.}\binits{Z.~D.}} \AND
  \bauthor{\bsnm{Silverstein},~\bfnm{J.~W.}\binits{J.~W.}}
(\byear{2004}).
\btitle{{CLT} for linear spectral statistics of large-dimensional sample
  covariance matrices}.
\bjournal{Ann. Probab.}
\bvolume{32}
\bpages{553--605}.
\end{barticle}
\endbibitem

\bibitem{bai2008central}
\begin{barticle}[author]
\bauthor{\bsnm{Bai},~\bfnm{Z.~D.}\binits{Z.~D.}} \AND
  \bauthor{\bsnm{Yao},~\bfnm{J.~F.}\binits{J.~F.}}
(\byear{2008}).
\btitle{Central limit theorems for eigenvalues in a spiked population model}.
\bjournal{Ann. Inst. Henri Poincar{\'e} Probab. Stat.}
\bvolume{44}
\bpages{447--474}.
\end{barticle}
\endbibitem

\bibitem{bai2012sample}
\begin{barticle}[author]
\bauthor{\bsnm{Bai},~\bfnm{Z.~D.}\binits{Z.~D.}} \AND
  \bauthor{\bsnm{Yao},~\bfnm{J.~F.}\binits{J.~F.}}
(\byear{2012}).
\btitle{On sample eigenvalues in a generalized spiked population}.
\bjournal{J. Multivariate Anal.}
\bvolume{106}
\bpages{167--177}.
\end{barticle}
\endbibitem

\bibitem{baik2005phase}
\begin{barticle}[author]
\bauthor{\bsnm{Baik},~\bfnm{J.}\binits{J.}},
  \bauthor{\bsnm{Ben~Arous},~\bfnm{G.}\binits{G.}} \AND
  \bauthor{\bsnm{P{\'e}ch{\'e}},~\bfnm{S.}\binits{S.}}
(\byear{2005}).
\btitle{Phase transition of the largest eigenvalue for nonnull complex sample
  covariance matrices}.
\bjournal{Ann. Probab.}
\bvolume{33}
\bpages{1643--1697}.
\end{barticle}
\endbibitem

\bibitem{2006Eigenvalues}
\begin{barticle}[author]
\bauthor{\bsnm{Baik},~\bfnm{J.}\binits{J.}} \AND
  \bauthor{\bsnm{Silverstein},~\bfnm{J.~W.}\binits{J.~W.}}
(\byear{2006}).
\btitle{Eigenvalues of large sample covariance matrices of spiked population
  models}.
\bjournal{J. Multivariate Anal.}
\bvolume{97}
\bpages{1382--1408}.
\end{barticle}
\endbibitem

\bibitem{Bao2022}
\begin{barticle}[author]
\bauthor{\bsnm{Bao},~\bfnm{Z.}\binits{Z.}},
  \bauthor{\bsnm{Ding},~\bfnm{X.}\binits{X.}},
  \bauthor{\bsnm{Wang},~\bfnm{J.}\binits{J.}} \AND
  \bauthor{\bsnm{Wang},~\bfnm{K.}\binits{K.}}
(\byear{2022}).
\btitle{Statistical inference for principal components of spiked covariance
  matrices}.
\bjournal{Ann. Statist.}
\bvolume{50}
\bpages{1144--1169}.
\end{barticle}
\endbibitem

\bibitem{bodnar2019testing}
\begin{barticle}[author]
\bauthor{\bsnm{Bodnar},~\bfnm{T.}\binits{T.}},
  \bauthor{\bsnm{Dette},~\bfnm{H.}\binits{H.}} \AND
  \bauthor{\bsnm{Parolya},~\bfnm{N.}\binits{N.}}
(\byear{2019}).
\btitle{Testing for independence of large dimensional vectors}.
\bjournal{Ann. Statist.}
\bvolume{47}
\bpages{2977--3008}.
\end{barticle}
\endbibitem

\bibitem{cheng2026supplement}
\begin{bmisc}[author]
\bauthor{\bsnm{Cheng},~\bfnm{H.}\binits{H.}},
  \bauthor{\bsnm{Li},~\bfnm{H.}\binits{H.}},
  \bauthor{\bsnm{Yin},~\bfnm{Y.}\binits{Y.}} \AND
  \bauthor{\bsnm{Zhang},~\bfnm{Z.}\binits{Z.}}
(\byear{2026}).
\btitle{Supplement to ``Limiting Eigen-Structure of Spiked Sample Covariance
  Matrices under Missing Observations''}.
\bnote{DOI to be provided by the typesetter}.
\end{bmisc}
\endbibitem

\bibitem{2018Generalized}
\begin{barticle}[author]
\bauthor{\bsnm{Jiang},~\bfnm{D.}\binits{D.}} \AND
  \bauthor{\bsnm{Bai},~\bfnm{Z.}\binits{Z.}}
(\byear{2021}).
\btitle{Generalized four moment theorem and an application to {CLT} for spiked
  eigenvalues of large-dimensional covariance matrices}.
\bjournal{Bernoulli}
\bvolume{27}
\bpages{274--294}.
\end{barticle}
\endbibitem

\bibitem{johnstone2001distribution}
\begin{barticle}[author]
\bauthor{\bsnm{Johnstone},~\bfnm{I.~M.}\binits{I.~M.}}
(\byear{2001}).
\btitle{On the distribution of the largest eigenvalue in principal components
  analysis}.
\bjournal{Ann. Statist.}
\bvolume{29}
\bpages{295--327}.
\end{barticle}
\endbibitem

\bibitem{2018Notes}
\begin{barticle}[author]
\bauthor{\bsnm{Johnstone},~\bfnm{I.~M.}\binits{I.~M.}} \AND
  \bauthor{\bsnm{Yang},~\bfnm{J.}\binits{J.}}
(\byear{2018}).
\btitle{Notes on asymptotics of sample eigenstructure for spiked covariance
  models with non-{G}aussian data}.
\bjournal{arXiv preprint arXiv:1810.10427}.
\end{barticle}
\endbibitem

\bibitem{li2024spectral}
\begin{barticle}[author]
\bauthor{\bsnm{Li},~\bfnm{H.}\binits{H.}},
  \bauthor{\bsnm{Pan},~\bfnm{G.}\binits{G.}},
  \bauthor{\bsnm{Yin},~\bfnm{Y.}\binits{Y.}} \AND
  \bauthor{\bsnm{Zhou},~\bfnm{W.}\binits{W.}}
(\byear{2024}).
\btitle{Spectral analysis of {G}ram matrices with missing at random
  observations: Convergence, central limit theorems, and applications in
  statistical inference}.
\bjournal{Ann. Statist.}
\bvolume{52}
\bpages{1254--1275}.
\end{barticle}
\endbibitem

\bibitem{Morales-Jimenez2021}
\begin{barticle}[author]
\bauthor{\bsnm{Morales-Jimenez},~\bfnm{D.}\binits{D.}},
  \bauthor{\bsnm{Johnstone},~\bfnm{I.~M.}\binits{I.~M.}},
  \bauthor{\bsnm{McKay},~\bfnm{M.~R.}\binits{M.~R.}} \AND
  \bauthor{\bsnm{Yang},~\bfnm{J.}\binits{J.}}
(\byear{2021}).
\btitle{Asymptotics of eigenstructure of sample correlation matrices for
  high-dimensional spiked models}.
\bjournal{Statist. Sinica}
\bvolume{31}
\bpages{571--601}.
\end{barticle}
\endbibitem

\bibitem{2007Asymptotics}
\begin{barticle}[author]
\bauthor{\bsnm{Paul},~\bfnm{D.}\binits{D.}}
(\byear{2007}).
\btitle{Asymptotics of sample eigenstructure for a large dimensional spiked
  covariance model}.
\bjournal{Statist. Sinica}
\bvolume{17}
\bpages{1617--1642}.
\end{barticle}
\endbibitem

\bibitem{waternaux1976asymptotic}
\begin{barticle}[author]
\bauthor{\bsnm{Waternaux},~\bfnm{C.~M.}\binits{C.~M.}}
(\byear{1976}).
\btitle{Asymptotic distribution of the sample roots for a nonnormal
  population}.
\bjournal{Biometrika}
\bvolume{63}
\bpages{639--645}.
\end{barticle}
\endbibitem

\bibitem{zhang2022asymptotic}
\begin{barticle}[author]
\bauthor{\bsnm{Zhang},~\bfnm{Z.}\binits{Z.}},
  \bauthor{\bsnm{Zheng},~\bfnm{S.}\binits{S.}},
  \bauthor{\bsnm{Pan},~\bfnm{G.}\binits{G.}} \AND
  \bauthor{\bsnm{Zhong},~\bfnm{P.~S.}\binits{P.~S.}}
(\byear{2022}).
\btitle{Asymptotic independence of spiked eigenvalues and linear spectral
  statistics for large sample covariance matrices}.
\bjournal{Ann. Statist.}
\bvolume{50}
\bpages{2205--2230}.
\end{barticle}
\endbibitem

\end{thebibliography}

\end{document}